\documentclass[12pt]{article}
\usepackage{latexsym}
\usepackage{epsfig}
\usepackage{enumerate}
\usepackage{amsmath,amsthm,amssymb}
\usepackage[active]{srcltx}
\usepackage[usenames]{color}
\usepackage{graphicx}
\definecolor{brown}{cmyk}{0, 0.72, 1, 0.45}
\definecolor{duck-shit}{cmyk}{0.3,0,1,0.4}
\definecolor{my-cyan}{cmyk}{1,0.4,0,0.2}
\definecolor{my-green}{cmyk}{.8,0,1,0}
\definecolor{grey}{gray}{0.5}
\def\red{\color{red} }

\def\sX{X}
\def\hZl{\wh{Z}}

\def\Zl{Z}
\def\Zu{Z}
\def\Tl{T^{low}}
\def\Tu{T^{up}}

\def\cF{{\cal F}}
\def\cG{{\cal G}}

\def\hd{\bar{d}}
\def\deg{\text{deg}}
\newcounter{rot}%\addtocounter{rot}{1}, \therot

\newcommand{\xbignote}[1]{}

\newcommand{\gap}[1]{\mbox{\hspace{#1 in}}}

\def\hp{\hat{\pi}}

\def\cN{{\cal N}}
\def\cL{{\cal L}}

\def\cLa{{\cL}_a}

\def\nn{\nonumber}
\def\a{\alpha} \def\b{\beta} \def\d{\delta} \def\D{\Delta}
\def\e{\epsilon}  \def\F{{\Phi}}  \def\g{\gamma}
  
\def\z{\zeta} \def\th{\theta}    \def\l{\lambda}
\def\La{\Lambda} \def\m{\mu} \def\n{\nu} \def\p{\pi}
\def\r{\rho}  \def\s{\sigma} 
\def\t{\tau} \def\om{\omega} \def\OM{\Omega}

\def\cD{{\cal D}}

\def\cC{{\cal C}}
\def\cB{{\cal B}}
\def\cT{{\cal T}}

\newtheorem{theorem}{Theorem}
\newtheorem{lemma}[theorem]{Lemma}
\newtheorem{corollary}[theorem]{Corollary}
\newtheorem{Remark}{Remark}

\newcommand{\proofstart}{{\bf Proof\hspace{2em}}}
\newcommand{\proofend}{\hspace*{\fill}\mbox{$\Box$}}
\newcommand{\W}[1]{W}

\def\cW{{\cal W}}
\def\cX{{\cal X}}
\newcommand{\ooi}{(1+o(1))}

\newcommand{\ul}[1]{\mbox{\boldmath$#1$}}

\newcommand{\wh}[1]{\widehat{#1}}

\newcommand{\rdup}[1]{{\left\lceil #1 \right\rceil }}
\newcommand{\rdown}[1]{{\left\lfloor #1\right \rfloor}}

\newcommand{\beq}[1]{\begin{equation}\label{#1}}
\newcommand{\eeq}{\end{equation}}
\newcommand{\blem}[1]{\begin{lemma}\label{#1}}
\newcommand{\elem}{\end{lemma}}
\newcommand{\bthm}[1]{\begin{theorem}\label{#1}}
\newcommand{\ethm}{\end{theorem}}

\newcommand{\brac}[1]{\left(#1\right)}

\newcommand{\sfrac}[2]{\frac{\scriptstyle #1}{\scriptstyle #2}}
\newcommand{\bfrac}[2]{\left(\frac{#1}{#2}\right)}
\def\half{\sfrac{1}{2}}
\def\cE{{\cal E}}
\newcommand{\rai}{\rightarrow \infty}
\newcommand{\ra}{\rightarrow}

\newcommand{\set}[1]{\left\{#1\right\}}
\def\sm{\setminus}
\def\seq{\subseteq}
\def\es{\emptyset}

\def\E{\mbox{{\bf E}}}
\def\Pr{\mbox{{\bf Pr}}}
\def\whp{{\bf whp}}
\def\Whp{{\bf Whp }}

\def\qs{{\bf qs}}

\newcommand{\ignore}[1]{}

\newcommand{\cA}{{\cal A}}

\def\var{{\bf Var}}
\begin{document}

\makeatletter
\title{Stationary distribution and cover time \\
of random walks on  random digraphs.
}

\author{Colin Cooper\thanks{Department of  Informatics,
King's College, University of London, London WC2R 2LS, UK}\and Alan
Frieze\thanks{Department of Mathematical Sciences, Carnegie Mellon
University, Pittsburgh PA15213, USA. Supported in part by NSF grant CCF0502793.}}

\maketitle \makeatother
\begin{abstract}
We study properties of a simple random walk on  the random digraph $D_{n,p}$ when
$np={d\log n},\; d>1$.

We prove that \whp\ the stationary
probability $\p_v$ of a vertex $v$ is asymptotic to $\deg^-(v)/m$
where $\deg^-(v)$ is the in-degree of $v$ and $m=n(n-1)p$ is the expected number of
edges of $D_{n,p}$.
If $d=d(n)\to\infty$ with $n$,  the stationary
distribution is asymptotically uniform \whp.

Using this result we prove that, for $d>1$,
\whp\ the cover time of $D_{n,p}$  is asymptotic to $d\log (d/(d-1))n\log n$.
If $d=d(n)\to\infty$ with $n$, then the cover time is asymptotic to $n\log n$.

\end{abstract}

\section{Introduction}
Let $D=(V,E)$ be a strongly connected digraph with $|V|=n$, and $|E|=m$.
For {the} simple random
walk $\cW_{v}=(\cW_v(t),\,t=0,1,\ldots)$ on $D$ starting at $v\in V$,
let $C_v$ be the expected time taken
 to visit every vertex of $D$.  The
{\em cover time} $C_D$ of $D$ is defined as $C_D=\max_{v\in V}C_v$.

For connected undirected graphs, the cover time
is well understood,
and has been extensively studied. It
is an old result of Aleliunas, Karp, Lipton, Lov\'asz and Rackoff
\cite{AKLLR} that $C_G \le 2m(n-1)$. It was shown by
 Feige \cite{Feige1}, \cite{Feige2},
that for any connected graph $G$, the cover time satisfies
$(1-o(1))n\log n\leq C_G\leq (1+o(1))\frac{4}{27}n^3$,
where $\log n$
is the natural logarithm.
An example of a graph achieving the lower
bound is the complete graph $K_n$  which has
cover time determined by the Coupon Collector problem.
The {\em lollipop} graph consisting of a path of length $n/3$ joined to a
clique of size $2n/3$ has cover time asymptotic to the upper bound of $(4/27)n^3$.

For  directed graphs  cover time is less well understood, and there are
strongly connected  digraphs with cover time exponential
in $n$.
An example of this is
 the digraph consisting of a directed cycle  $(1,2,...,n,1)$,
and edges $(j,1)$, from vertices  $j=2,...,n-1$.
Starting from vertex $1$, the expected time for a random walk to
reach vertex $n$ is {$\Omega(2^{n})$}.

In earlier papers, we investigated the cover time of various classes
of {(undirected)} random graphs, and derived precise results for their
cover times. The main results  can be summarized as
follows:
\begin{itemize}
\item \cite{CFC} If $p=d\log n/n$ and $d>1$ then \whp\ $C_{G_{n,p}}\sim
  d\log\bfrac{d}{d-1}n\log n$.
\item \cite{CFgiant,CFgiant1} Let $d>1$ and let $x$ denote the solution in
  $(0,1)$ of $x=1-e^{-dx}$.
Let $X_g$ be the giant component of $G_{n,p},\,p=d/n$. Then \whp\
$C_{X_g}\sim \frac{dx(2-x)}{4(dx-\log d)} n (\log n)^2$.
\item \cite{CFreg} If {$r\geq 3$ is a constant and }$G_{n,r}$
denotes a random $r$-regular graph on vertex set
  $[n]$ with $r\geq3$ then \whp\ $C_{G_{n,r}}\sim \frac{r-1}{r-2}n\log n$.
\item \cite{CFweb} If {$m\geq 2$ is constant and }$G_m$ denotes a
{\em preferential attachment graph} of
  average degree  $2m$ then \whp\ $C_{G_m}\sim \frac{2m}{m-1}n\log n$.
\item \cite{CFgeom} If $k\geq 3$ and $G_{r,k}$ is a random geometric graph in $\Re^k$ of
ball size $r$ such that the expected degree of a vertex is asymptotic to $d\log n$, then \whp\
 $C_{G_{r,k}}\sim   d\log\bfrac{d}{d-1}n\log n$.
\end{itemize}
A few remarks on notation:
We use the notation $a(n) \sim b(n)$ to mean that $a(n)/b(n) \ra 1$ as $n \rai$.
Some inequalities
in this paper only hold for large $n$. We  assume henceforth that $n$ is sufficiently large for
all claimed inequalities to hold.
All \whp\ statements in this paper are relative to the  class of random digraphs $D_{n,p}$
under discussion, and not the random walk.

In this paper we turn our attention to the cover time of random directed graphs.
Let $D_{n,p}$ be the random digraph with vertex set $V=[n]$ where
each possible directed edge $(i,j),\,i\neq j$ is {independently} included with
probability $p$. It is known that if $np=d\log n=\log n+\g$ where
$\g=(d-1)\log n\to \infty$ then $D_{n,p}$ is strongly connected
\whp. If $\g$ as defined tends to $-\infty$ then \whp\ $D_{n,p}$ is
not strongly connected. As we do not have a direct reference to this
result, we next give a brief proof of this.
 It is easy to show  that if $np=\log n-\g$ where $\g\to\infty$,
 there are vertices of in-degree zero \whp. On the other hand,
if $np=\log n+\g$ where $\g\to\infty$ then \cite{dham} shows that
the random digraph is Hamiltonian and hence strongly connected.
Strong connectivity for $np=\log n+\g$ where $\g\to\infty$ also follows directly from
the proof of \eqref{lowpi}.

We determine the cover time of $D_{n,p}$ for  values of $p$ at or above the
threshold for strong connectivity.
\begin{theorem}\label{thmain}\
Let $np=d\log n$ where {$d=d(n)$ is such that $\g=np-\log n\to\infty$}.
Then \whp
$$C_{D_{n,p}}\sim d\log\bfrac{d}{d-1}n\log n.$$
\end{theorem}
%{Here $X=X(n)\sim Y=Y(n)$ if $\lim_{n\to\infty}X(n)/Y(n)=1$.}

Note that if $d=d(n)\to\infty$ with $n$, then we have $C_{D_{n,p}}\sim n\log n$.

The method  we use to find the cover time of $D_{n,p}$
requires us to know the stationary distribution
of the random walk.
For an undirected
graph $G$,  the stationary distribution is
$\p_v=\deg(v)/2m$, where $\deg(v)$ denotes {the} degree of vertex $v$, and $m$ is the
number of edges in $G$. For  a digraph $D$,
 let $\deg^-(v)$ denote the in-degree of vertex $v$, $\deg^+(v)$ denote the out-degree,
 and let  $m$ be the
number of edges in $D$. For  strongly connected
  digraphs in which  each vertex $v$
has  in-degree equal to out-degree ($\deg^-(v)=\deg^+(v)$), then
$\pi_v=\deg^-(v)/m$. For general digraphs, however, there is no
simple formula for the stationary distribution.
Indeed, there may not be a unique stationary
measure.
The main technical task of this paper is to find
good estimates for $\pi_v$ in the case of $D_{n,p}$. Along the way,
this implies uniqueness of the stationary measure \whp.

We summarize our
result concerning the stationary distribution
in Theorem \ref{lemsteady} below.
For a given vertex $v$, define a quantity $\varsigma^*(v)$, which {in essence depends on}
 the in-neighbour $w$ of $v$ with minimum out-degree:
\beq{varsig}
\varsigma^*(v)=\max_{w \in N^-(v)} \set{\frac{\deg^-(w)}{\deg^+(w)}}.
\eeq
\begin{theorem}\label{lemsteady}\
Let $np=d\log n$ where $d=d(n)$ is such that $\g=np-\log n\to\infty$.
Let $m=n(n-1)p$.
Then \whp, for all $v \in V$,
$$\p_v\sim \frac{\deg^-(v)+ \varsigma^*(v)}{m}.$$
\end{theorem}
We note the following special cases.
\begin{Remark}\label{rem1}
We prove in Lemma \ref{newlem0} that \whp\
$\varsigma^*(v)=o(\deg^-(v))$ for almost all vertices $v$.
For these vertices, the $\varsigma^*(v)$ term can be absorbed into the
 error term of $\pi_v$.
\end{Remark}
\begin{Remark} If $\g=\om(\log\log n)$ then \whp\  $\varsigma^*(v)=o(\deg^-(v))$
for   all vertices $v$.
In particular when $d=1+\d, \; \d >0$ and constant then
the minimum out-degree  is $\Omega(\log n)$,
in which case, $\p_v \sim {\deg^-(v)}/{m}$.
\end{Remark}
\begin{Remark}
If $d=d(n)\to\infty$ with $n$, \whp\  the stationary
distribution of $D_{n,p}$ is $\pi_v \sim 1/n$.
\end{Remark}

\section{Outline of the paper}
At the heart of our approach to the cover time is the following claim: Suppose that
$T$ is a "mixing time" for a simple random walk and $\ul A_v(t)$ is the event that $\cW_u$
does not visit $v$ in steps $T,T+1,\ldots,t$. Then, essentially,
\beq{heart}
\Pr(\ul A_v(t))\sim e^{-t\p_v/R_v}.
\eeq
Here $R_v\geq 1,\,v\in V$ is the expected number of
visits/returns to $v$ by the walk $\cW_v$ within $T$ time steps. This is the content of
Lemma \ref{MainLemma} and we have used it to prove our previous results on this topic.
Given \eqref{heart} we can
estimate the cover time from above via
$$C_u \le t+1+\sum_v\sum_{s \ge t}\Pr(\ul A_v(s)).$$
This is \eqref{shed1} and we have used this inequality previously. Here $C_u$ is the
expected time for $\cW_u$ to visit every vertex. It is valid for arbitrary $t$ and we get
our upper bound for $C_D$ by choosing $t$ large enough so that the double sum is $o(t)$.

We estimate the cover time from below by using the Chebyshev inequality.
We choose a set of vertices $V^{**}$ that are candidates for taking
a long time to visit and estimate the
expected size of the set $V^\dagger$ of vertices in $V^{**}$
that have not been visited within our estimate of the cover time.
We show that $\E |V^\dagger|\to\infty$.
We then take pairs of vertices $v,w\in V^{**}$ and contract them to a single vertex $\g$
and then use \eqref{heart} to show that $\Pr(\ul A_\g(t))\sim \Pr(\ul A_v(t))\Pr(\ul A_w(t))$.

The main problem here is that we do not know $\p_v$ and much of the paper is devoted to proving that,
essentially, \whp,
\beq{steady}
\p_v\sim \frac{\deg^-(v)}{m}\qquad\text{for all }v\in V.
\eeq
Our proof of this leads easily to a claim that \whp\ $T=O(\log^2n)$ and we will find then find that
it is easy to prove that $R_v=1+o(1)$ for all $v\in V$.

We approximate the stationary distribution
$\ul \pi$ using the expression
 $\ul \pi = \ul \pi P^k$, where $P$ is the transition matrix.
For suitable choices of $k$ we find we can bound
$$P_x^{(k)}(y)=\Pr(\cW_x(k)=y)$$
from above and below by values independent of $x$
and obtain, essentially,
$$P_x^{(k)}(y)\sim \frac{\deg^-(y)}{m}$$
an expression independent of $x$. \eqref{steady} follows easily from this.

To estimate $P_x^{(k)}(y)$ from below we proceed as follows: We let $k=2\ell=\frac23\log_{np}n$.
We consider two Breadth First Search trees of depth $\ell$. $T^{low}_x$ branches out from
$x$ to depth $\ell$ and $T^{low}_y$ branches into $y$ from depth $\ell$.
Almost all of the walk measure associated with walks
of length $2\ell+1$ from $x$ to $y$ will go from $x$ level by level
to the boundary of $T^{low}_x$, jump across
to the boundary of $T^{low}_y$ and then go level by level to $y$. We analyse such walks
and produce a lower bound.

To estimate $P_x^{(k)}(y)$ from above we change the depths of the out-tree from $x$
and the in-tree to $y$. This eliminates some complexities. In computing the lower bound,
we ignored some paths that take more circuitous routes from $x$ to $y$ and we have to
show that these do not add much in walk measure.

The structure of the paper is now as follows:
Section \ref{mainlem} describes Lemma \ref{MainLemma} that we have often used
before in the analysis
of the cover time.
%Section \ref{useful} proves a technical lemma on a branching process,
%Lemma \ref{LemWXY}, that will be used to
%analyse the growth of Breadth First Search trees $T^{low}_y,T^{up}_y$.
Section \ref{structure} establishes many structural properties of
$D_{n,p}$.
In Section \ref{SteadyState} we prove the lower and upper
bounds given in Theorem \ref{lemsteady}.
These bounds hold for any digraph
with the high probability structures elicited in Section \ref{structure}.
Sections \ref{structure} and \ref{SteadyState}, which form the main body of this paper,
are first proved under the assumption that $2 \le d \le n^{\d}$, for some small $\d>0$,
an assumption we refer to as {\bf Assumption 1}.
In Section \ref{remove}, we extend the proof of Theorem \ref{lemsteady} by removing Assumption 1.
Section \ref{smallsec} is short and establishes that the
conditions of Lemma \ref{MainLemma} hold. To do this, we use a bound on the mixing time,
 based on results obtained in Sections \ref{SteadyState}, \ref{remove}.
Finally, in Section \ref{cta} we establish the \whp\ cover time, as given in Theorem \ref{thmain}.

\section{Main Lemma}\label{mainlem}
In this section $D$ denotes a fixed strongly connected digraph with $n$ vertices.
A random walk $\cW_{u}$  is started from a vertex $u$.
Let $\cW_{u}(t)$ be the vertex
reached at step $t$, let $P$ be the matrix of transition probabilities of the walk and let
$P_{u}^{(t)}(v)=\Pr(\cW_{u}(t)=v)$. We assume that the random walk $\cW_{u}$ on $D$ is ergodic with
 stationary distribution $\pi$.

Let
$$d(t)=\max_{u,x\in V}|P_{u}^{(t)}(x)-\pi_x|,$$
and
let $T$ be a positive integer such that for $t\geq T$
\begin{equation}\label{4}
\max_{u,x\in V}|P_{u}^{(t)}(x)-\pi_x| \leq n^{-3} .
\end{equation}
Consider the  walk $\cW_v$, starting
at {vertex} $v$. Let $r_t={r_t(v)=}\Pr(\cW_v(t)=v)$ be the probability  that this  walk
returns to $v$ at step $t = 0,1,...$\ .
Let
\begin{equation}
\label{Qs}
R_T(z)=\sum_{j=0}^{T-1} r_jz^j
 \end{equation}
and
let
$$R_v=R_T(1).$$

\begin{lemma}\label{MainLemma}
\label{L3}
Fix a vertex $u\in V$ and for $v\in V$ and $t\geq T$ let $\ul A_v(t)$ be the event that $\cW_u$
does not visit $v$ in steps $T,T+1,\ldots,t$.
Suppose that
\begin{description}
\item[(a)]
For some constant $\th >0$, we have
$$\min_{|z|\leq 1+\l}|R_T(z)|\geq \th.$$
\item[(b)]$T^2\pi_v=o(1)$ and $T\pi_v=\Omega(n^{-2})$ for all $v\in V$.
\end{description}
Let $K>0$ be a sufficiently large absolute constant and let
\begin{equation}
\label{lamby}
\l=\frac{1}{KT}.
\end{equation}
Then, with
\begin{equation}\label{pv}
p_v=\frac{\pi_v}{{R_v}(1+O(T\pi_v))},
\end{equation}
we have that
for all $v\in V$ and $t\geq T$,
\begin{equation}
\label{frat}
\Pr(\ul A_v(t))=\frac{(1+O(T\p_v))}{(1+p_v)^{t}} +O(T^2 \p_v e^{-\l t/2}).
\end{equation}
\end{lemma}

%%%end of ignore

\section{Structural Properties of $D_{n,p}$}\label{structure}
In this section we gather together some \whp\ properties of $D_{n,p}$
needed for the proof of Theorem \ref{lemsteady}.
Some are quite elaborate and so we will try to motivate them where we can.

We stress that throughout this section, the probability space is the space of $D_{n,p}$
and not the space of walks on an instance.

Once we complete this section however, we can concentrate on estimating the cover time
of a digraph with the given properties, as long as Assumption 1 of \eqref{Ass1} holds.
For large $p$ outside Assumption 1, the proof is quite simple.

\subsection{Degree Sequence etc}\label{secdegree}

{\bf Chernoff Bounds}
The following inequalities are used extensively throughout this paper.
Let   $Z=Z_1+Z_2+\cdots Z_N$ be the sum of  the independent random variables
$0\leq Z_i\leq 1,\,i=1,2,\ldots,N$ with $\E(Z_1+Z_2+\cdots+Z_N)=N\m$. Then
\begin{eqnarray}
\Pr(|Z-N\m|\geq \e N\m)&\leq&2e^{-\e^2N\m/3}.\label{Cher1}\\
\Pr(Z\geq \a N\m)&\leq&(e/\a)^{\a N\m}.\label{Cher2}
\end{eqnarray}
{For proofs see for example Alon and Spencer \cite{AS}.}

The next lemma gives some properties of  the degree sequence of $D_{n,p}$.
The lemma
can  be proved by the use of the first and second moment
methods (see \cite{CFC} for very similar calculations).

Let $np=d \log n$ and let
\beq{D0}
\D_0=C_0 np\ where\ C_0=30.
\eeq
\begin{lemma}\label{vertexdegrees}\
\begin{enumerate}[(i)]
\item First assume that $np=d\log n$ where $1<d=O(1)$ and $(d-1)\log n\to\infty$.
Let
$$\overline{D}(k)=n\binom{n-1}{k}p^k(1-p)^{n-1-k}$$
denote the expected number of vertices $v$ with $\deg^-(v)=k\leq\D_0$.
{Note that
\beq{aproxval}
\overline{D}(k)\leq \frac{2}{n^{d-1}}\bfrac{ne{p}}{k}^k.
\eeq}
Let $D(k)$ denote the actual
number {of vertices of in-degree $k$}, and let
\begin{eqnarray*}
K_0&=&\{k\in [1,\D_0]:\;\overline{D}(k)\leq (\log n)^{-2}\}.\\
K_1&=&\{1\leq k\leq 15:\;(\log n)^{-2}\leq
\overline{D}(k)\leq \log\log n\}.\\
K_2&=&\{k\in [16,\D_0]:\;(\log n)^{-2}\leq \overline{D}(k)\leq (\log
n)^2\}.\\
K_3&=&[1,\D_0]\setminus (K_0\cup K_1\cup K_2).
\end{eqnarray*}
The degree sequence has the following properties.
\begin{description}
\item[(a)] If $d-1\ge (\log n)^{-1/3}$ then
$$K_1=\emptyset, \qquad \min \{k \in K_2\} \ge (\log n)^{1/2},\qquad |K_2|= O(\log \log n).$$
\item[(b)] {The following holds \whp:} For all degrees $k \in [1,\D_0]$.
\begin{eqnarray}
k\in K_0, & & D(k)=0,\nonumber\\
k\in K_1, &  & D(k) \leq (\log\log n)^2,\label{K1val}\\
k\in K_2, &  & D(k) \leq (\log n)^4,\label{K2val}\\
k\in K_3, &&
\frac{\overline{D}(k)}{2}\leq D(k)\leq
2\overline{D}(k).\label{K3val}
\end{eqnarray}
\end{description}
\item Suppose that $1<d\leq n^\d$ where $\d$ is a small positive constant.
Let $k^*=\rdup{(d-1)\log n}$.
Let
$$V^*=\set{v\in V:\;\deg^-(v)=k^*\text{ and }\deg^+(v)=
k^\dagger=\rdup{d\log n}}$$
and let $\g_d=(d-1)\log\bfrac{d}{d-1}$. Then \whp
$$|V^*|\geq
\frac{n^{\g_d}}{10d\log n}.$$
\item
Let $\cD$ be the event
\beq{calD}
\set{\exists v\in V:\;\deg^+(v)\geq\D_0\ or\ \deg^-(v)\geq\D_0},
\eeq
then
\beq{largedegree}
\Pr(\cD)\leq n^{-10}e^{-10np}.
\eeq
\item The number of edges $|E(D_{n,p}|\sim m=n(n-1)p$ \whp.
\item $\deg^\pm(v)\sim np$ for all $v\in V$ \whp\ if $d\to\infty$.
\end{enumerate}
\end{lemma}
\proofstart
We will only give an outline proof of (ii)
as the other claims have (essentially) been proved in \cite{CFC}.
We have
$$\E(|V^*|)=n\binom{n-1}{k^*}\binom{n-1}{k^\dagger}p^{k^*+k^\dagger}(1-p)^{2n-2-k^*-k^\dagger}$$
Using, $\binom{n}{k}\geq \frac{1}{3k^{1/2}}\bfrac{ne}{k}^k$, which is obtainable from Stirling's
approximation, we see that
\beq{EV*}
\E(|V^*|)\geq (1-o(1))\frac{1}{9(d-1)^{1/2}d^{1/2}\log n}\bfrac{d}{d-1}^{(d-1)\log n}.
\eeq
The Chebyshev inequality will show that $|V^*|$ is concentrated around its mean, provided we verify
that the mean tends to $\infty$.

If $d>2$ then
$$\bfrac{d}{d-1}^{(d-1)\log n}\geq
\exp\set{(d-1)\log n\brac{\frac{1}{d-1}-\frac{1}{2(d-1)^2}}}\geq n^{1/2}.$$
Now $d\leq n^\d$ and so $\E(|V^*|)\to\infty$ follows from \eqref{EV*}. If $d=1+\e\leq 2$ and
$\e$ is bounded away from zero then so is $\bfrac{d}{d-1}^{d-1}$ and so $\E(|V^*|)=n^{\Omega(1)}$.
So now suppose that $\e=\frac{\om}{\log n}$ where $\om=\om(n)\to\infty$ and $\om=o(\log n)$. Then,
$\bfrac{d}{d-1}^{d-1}\geq \bfrac{\log n}{\om}^\om$. If $\om\geq \log^{1/2}n$ then
$\bfrac{d}{d-1}^{d-1}\geq e^{\log^{1/2}n}$ and if $\om\leq \log^{1/2}n$ then
$\bfrac{d}{d-1}^{d-1}\geq \log^{\om/2}n$. In either case
\beq{EVL}
\E(|V^*|)\geq \log^{\th}n\text{ where $\th=\th(n)\to\infty$}
\eeq
\proofend
\subsubsection{Assumption 1: a convenient restriction}

We will first carry out the
main body of the proof under the  following assumption:
\beq{Ass1}
\gap{-2}{\bf Assumption\ 1}:\gap{1} 2\leq d\leq n^\d.
\eeq
Here $0<\d\ll 1$ is some small fixed positive constant.
We note that our choice of the value $d\ge 2$ is somewhat arbitrary,
and any constant larger than 1 would suffice.
We wait until Section \ref{remove} to remove Assumption 1.
The proof for $d>n^\d$ is much
simpler and is given separately in Section \ref{bigdegree}.
The proof for $1 <d \le 2$ is given in Section \ref{Smavdeg}.

\ignore{
Suppose now that we make the following assumption:\\
{\bf Assumption 1}: $d-1=\Omega(1)$ i.e. $d$ is at least a positive
constant strictly larger than one.

Let $c_b=d'/2$ and $C_b=d''+1$ where $d'<d''$
are the roots of $x\log(e/x)=1-{\log n}/{np}$. Let $I_b=[c_bnp,C_bnp]$. Then, with $I=I_b$,
$$\Pr\brac{\bigcap_{v\in V}\cD^+(v)\cap \bigcap_{v\in V}\cD^-(v)}=\r_b=1-O(n^{-\psi})$$
where $\psi=\psi(d)>0$.
}

Under Assumption 1 and $d=O(1)$ there is a constant $c>0$ and an interval
\beq{eqII}
I=[c_0np,\D_0],
\eeq
such that if {$\n\in [3n/4,n]$} then there exists $\g=\g(c)>0$ such that
\beq{eqI}
\Pr(Bin(\n,p)\in I)=1-o(n^{-1-\g}).
\eeq
When $d\to\infty$ we can take $c=0.999$ and $C_0=1.001$.

Let $\cE^+_S$ (resp. $\cE^-_S$) be the event that the in-degree (resp. out-degree) of
all vertices in $S\subseteq V$ are in the interval $I$. Thus e.g.
\beq{EsetS}
\cE^+_S= \set{D_{n,p}: \forall v \in S\subseteq V, \deg^+(v) \in [c_0np,\D_0]}.
\eeq
Let $\cE_S=\cE^+_S\cap \cE^-_S$. Then
for any $S\subseteq V$ we have
\beq{cE}
\Pr(\cE_S)=1-O(n^{-\g}).
\eeq

%\subsection{Breadth First Search Trees}

\subsection{Properties needed for a lower bound on the stationary distribution}\label{lbd}

The calculations in this section are made under Assumption 1.

Fix vertices $x,y$ where $x=y$ is allowed.
Most  short random  walks from  vertex $x$ to  vertex $y$  take the
form of a simple directed path, or cycle if $x=y$. We can count such paths
(or cycles) with the help of  a breadth first out-tree $\Tl_x$ rooted at $x$,
and a breadth first in-tree  $\Tl_y$ rooted at $y$.
We build these trees to depth $\ell$, where
\begin{equation}\label{ell}
\ell=\rdown{\frac23\log_{np}n}.
\end{equation}

For a vertex $v$ let $N^-(v)$ be the set of in-neighbours of $v$ and
for a set $S$, let $N^-(S)=\bigcup_{v\in S}N^-(v)$.
Define $N^+(v),N^+(S)$  similarly with respect to out-neighbours.

{\bf Construction of  in-tree $\Tl_y$.}

For fixed $y\in V$,
we build a tree $T_y=\Tl_y$ rooted at $y$, in
breadth-first fashion. Define
$Y_0=\set{y}$.
The tree
$\Tl_y$ has level sets $ Y_i, i=0,1,...,\ell$, and vertex set $Y=\cup_{i=0}^{\ell} Y_i$.
Let $Y_{\le i} = \cup_{j \le i} Y_j$.

Let $T_y(i)$ be the tree consisting of the first $i$
levels of the breadth first tree $T_y(\ell)=\Tl_y$.
Given $T_y(i)$ we construct $T_y(i+1)$ by adding the in-neighbours of $Y_i$ in $V \sm Y_{\le i}$.
To remove ambiguity, the vertices of $Y_i$ are processed in increasing order
of vertex label. Let this order be $(v_1,v_2,...,v_{|Y_i|})$.
For $v \in Y_i$, the set $N_T^-(v)$ is the subset of $Y_{i+1}\cap N^-(v)$
whose edges in the tree $T_y(i+1)$ point to $v$.
 Formally these sets are defined as follows:
$N_T^-(v_1)=N^-(v_1) \sm Y_{\le i}$, and in general
$N_T^-(v_k)= N^-(v_k) \sm (Y_{\le i} \; \cup N^-(v_1,...,v_{k-1}))$.

Thus if $v \in Y_i$ and $w\in Y_{i+1}$ and $(w,v)$ is an edge of $T_y(i+1)$,
then $v=v_k$ is the first out-neighbour of $w$ in $Y_i$ in the order $(v_1,...,v_k)$.
As $w \in Y_{i+1}$, there are no edges from $w$ to $Y_{\le i-1}$, and thus
no edges between  $w$ and  $Y_{\le i-1} \cup \{v_1,...,v_{k-1}\}$.

Let $\deg_T^-(v)$ denote the
in-degree of $v\in Y$ in $\Tl_y$.
If $\deg_T^-(v)>0$ for all $v \in Y_{\le \ell-1}$, we say {\em the construction of $\Tl_y$ succeeds}.

Associated with this construction  of  $\Tl_y$ is a set of
 parameters and random variables.

\begin{itemize}
\item For $v_j \in Y_i$, let $\s(v_j) = |V \sm [Y_{\le i} \; \cup N_T^-(v_1,...,v_{j-1})]|$.
Thus  $\s(v_j)$ is the number of vertices not in $T_y$
after all in-neighbours of $v_1,...,v_{j-1}$ have been added to $T_y$.
\item Let $B(v_j) \sim Bin(\s(v_j),p)$. Thus
$B(v_j)=|N_T^-(v_j)|=\deg_T^-(v_j)$, the in-degree of $v_j$ in $T_y$.
\item Let $\s'(v_j)=|V \sm [Y_{\leq i-1} \cup  \{v_1,...,v_{j-1}\}|$.
\item Let $D \sim 1+Bin(\s'(v_j),p)$, and let $D(j,k),\; k=1,...,B(v_j)$ be independent copies of $D$.
\end{itemize}
The interpretation of the random variable $D(j,k)$ is as follows.
If $w_k \in N_T^-(v_j)$ then $D(j,k)$ is the out-degree of $w_k$ in $D_{n,p}$.
The one arises from $(w_k,v_j)$
being the first edge from $w_k$ to $Y_i$.

{\bf Construction of out-tree $\Tl_x$.}
Given the set of vertices $Y$ of $\Tl_y$, we define
$X_0=\set{x},X_1,\ldots,X_\ell$ where $X_{i+1}=N^+(X_i)\sm
(Y\cup X_0\cup\cdots\cup X_i)$ for $0\leq i<\ell$.
If $w\in X_{i+1}$ is the out-neighbour of
more than one vertex of $X_i$, we only keep the edge $(z,w)$
with $z$ as small as possible as in the construction of $\Tl_y$.
Let $X=\bigcup_{i=0}^\ell X_i$ and
let $\Tl_x$ denote the BFS tree
constructed in this manner.
Let $\deg_T^+(v)=|N^+_T(v)|$ denote the
out-degree of $v\in X$ in $\Tl_x$.
Similarly to the construction of $\Tl_y$, the value of $\deg^+_T(v)$ is
given by a random variable $B(v) \sim Bin(\s(v),p)$.
If $\deg_X^+(v)>0$ for all $v \in X_{\le \ell-1}$,
we say the construction of $\Tl_x$ succeeds.

We gather together a few facts about $\Tl_x,\Tl_y$ that we need for the proofs of this section.

\begin{lemma}\label{treelemma}\
\begin{description}
\item[(i)]With probability $1-O(n^{-\g})$, the construction of $\Tl_x,\Tl_y$ succeeds
for all $x,y \in V$.
\item[(ii)]With probability $1-O(n^{-\g})$,\\
For all $x$, and for all $v \in X_{\le \ell-1}$, $\deg^+_T(v) \in [c_0np(1-o(1)), C_0np]$,
\\
For all $y$ and for all $v \in Y_{\le \ell-1}$, $\deg^-_T(v) \in [c_0np(1-o(1)), C_0np]$.
\item [(iii)] Given $\cE^+_x,\cE^-_y$, then  for $i \le \ell$,
$|X_i| \sim \deg^+(x) (np)^{i-1}, \;|Y_i| \sim \deg^-(y) (np)^{i-1}$
\qs\footnote{A sequence of events $\cA_n,n\geq 0$ hold \lq quite surely\rq (\qs) if
$\Pr(\cA_n)=1-O(n^{-K})$ for any constant $K>0$.}
\end{description}
\end{lemma}
\proofstart
We give  proofs for $\Tl_x$, the proofs for $\Tl_y$ are similar.

{\bf Part (i), (ii).}
Let $X=\set{x_0=x,x_1,\ldots,x_N}$ where $x_i$ is the $i$-th vertex added to $\Tl_x$.
For $x_j\in X$, let $f(x_j)=|N^+(x_j)\cap(Y\cup \set{x_0,x_1,\ldots,x_{j-1}})|$.
Thus $\deg^+(v)=\deg^+_T(v)+f(v)$.

We can bound {$f(v)$ by the binomial $Bin(N_X,p)$ } where
$N_X=|Y|+|X|$. This is true even after constructing $\Tl_x,\Tl_y$,
because the out-edges of $v$ counted by $f(v)$
have not been exposed.
Assuming $\neg\cD$, see \eqref{largedegree}, we have
$$N_X \leq 2\sum_{i=1}^{\ell}\D_0^i=n^{2/3+o(1)}.$$
Using the Chernoff bound \eqref{Cher2}, we have {with $\om=\log^{1/2}n$} that
\beq{just1}
\Pr\brac{ f(v) \ge \frac{np}{\om}} \le
\Pr\brac{Bin(n^{2/3+o(1)},p)\ge \frac{np}{\om}}
=O(n^{-10}).
\eeq
The event $\bigcup_{x \in v}\cE^+_{X_{\le \ell-1}} \seq\cE^+_V$
and the latter holds with probability $1-O(n^{-\g})$.
Thus given \eqref{just1},  and $\cE^+_{V}$ we have
$\deg^+_T(v)>0$ for $v\in X_{\le \ell-1}$ for all $x \in V$.
In summary, \whp\ the construction of
$\Tl_x$ succeeds for all $x \in V$, and
 $\deg^+_T(v) \in [c_0np(1-o(1)), C_0np]$ for all $v \in X_{\le \ell-1}$
in all trees $\Tl_x$, $x \in V$.

{\bf Part (iii).}
By construction $\Tl_x$ was made after $\Tl_y$, so $|X_i|$ depends on $\Tl_x$ and $\Tl_y$.
Assume $\cE^+_x = \set{\deg^+(x) \in  I=[c_0np,C_0np]}$, and  that $|Y| \le n^{2/3+o(1)}$.
(Strictly speaking we should verify that $|Y| \le n^{2/3+o(1)}$
before considering $\Tl_x$. On the other hand,
the proof we present now can be applied to $\Tl_y$).
For $i\geq 1$, $|X_{i+1}|$ is distributed as $Bin(n-o(n),1-(1-p)^{|X_i|})$.
The number of trials $n-o(n)$ is based on the inductive assumption that
$|X_{j+1}|=(1+o(1))n|X_j|p$ and that $|X_{\leq j}|p=o(1)$. That these assumptions hold
\qs.
follows from the Chernoff bounds. We thus have that
\beq{xx1}
|X_\ell|\sim \deg^+(x)(np)^{\ell-1},\qquad\qs.
\eeq
\proofend

For $u\in X_i$ let $P_u$ denote the path of length $i$ from $x$ to $u$ in $\Tl_x$ and
$$\a_{i,u}=\prod_{\substack{w\in P_{u}\\w \ne u}}\frac{1}{\deg^+(w)}.$$
In the  event that  the construction of $\Tl_x$ fails to complete to depth $\ell$,
let $\sum_{u\in X_\ell}\a_{\ell,u}=0$.

Similarly, for $v\in Y_i$  let $Q_{v}$ denote the path from $v$ to $y$ in $\Tl_y$ and
\begin{equation}\label{betaiv}
\b_{i,v}=\prod_{\substack{w\in Q_{v}\\ w \ne y}}\frac{1}{\deg^+(w)}.
\end{equation}
In the  event that the construction of $\Tl_y$ fails to complete to depth $\ell$,
let $\sum_{v\in Y_\ell}\b_{\ell,v}=0$.

Let
\begin{equation}\label{Given}
\Zl(x,y)=Z^{low}(x,y)=\sum_{\substack{u\in X_\ell\\v\in
    Y_\ell}}\a_{\ell,u}\b_{\ell,v}\frac{1_{uv}}{\deg^+(u)}
\end{equation}
where $1_{uv}$ is the indicator for the existence of the edge $(u,v)$
and we take $\frac{1_{uv}}{\deg^+(u)}=0$ if $\deg^+(u)=0$.
Note that $Z(x,y)=0$ if we fail to construct $\Tl_x$ or $\Tl_y$.

\begin{Remark}
The importance of the quantity $\Zl(x,y)$ lies in the fact that it is
a lower bound on the probability that $\cW_x(2\ell+1)=y$.

The aim of the next few lemmas is to prove (see Lemma \ref{xlem2}), under Assumption 1, that with $m=n(n-1) p$,
\beq{promise}
\Pr\brac{\exists_{x,y \in V} \text{ such that } Z(x,y)
\ne \frac{(1+o(1))\deg^-(y)}{m} }=O(n^{-\g}).
\eeq
The first two lemmas give \whp\  bounds for $\sum_{u \in X_\ell}
\a_{\ell,u}$, $\sum_{v \in Y_\ell}\b_{\ell,v}$ respectively, as used in
 the third lemma and its corollary.
\end{Remark}

\begin{lemma}\label{xlem1}
Let
\beq{cA1}
\cA_1(x,y)=\set{1-\e_X \le \sum_{u\in X_\ell}\a_{\ell,u} \le 1}.
\eeq
For $\e_X=O(1/\sqrt{\log n})$,
\begin{equation}
\Pr\brac{\neg\cA_1(x,y)\cap \cE^+_x}=o(n^{-10}).\label{LSa}
\end{equation}

\end{lemma}
\proofstart

For $u\in X_\ell$ let $xP_u=(u_0=x,u_1,\ldots,u_\ell=u)$
denote the path from $x$ to $u$ in $\Tl_x$.
For the random walk on the digraph  $\Tl_x$, starting at $x$; $X_{\ell}$
is reached with probability $\F=1$ in exactly $\ell$ steps,
after which the walk halts. Thus
\beq{mn}
1=\F= \sum_{u\in X_\ell}\prod_{\substack{v\in P_u\\v \ne u}}\frac{1}{\deg^+_T(v)}
\geq \sum_{u\in X_\ell}\a_{\ell,u}.
\eeq
We assume that the construction of $\Tl_x$ succeeds, and
that $\deg^+_T(v)>0$ for $v\in X_{\le \ell-1}$, as established in Lemma \ref{treelemma}.
In the notation of that lemma, $\deg^+(v)=\deg^+_T(v)+f(v)$.
Now
\begin{eqnarray*}
\F&=&\sum_{u\in X_\ell}\prod_{v\in P_u}\frac{1}{\deg^+(v)-f(v)}\\
&=&\sum_{u\in X_\ell}\brac{\prod_{v\in P_u}\frac{1}{\deg^+(v)}}
\brac{\prod_{v\in P_u}\frac{1}{1-f(v)/\deg^+(v)}}\\
&=&\sum_{u\in X_\ell}\a_{\ell,u}\brac{\prod_{v\in P_u}\frac{1}{1-f(v)/\deg^+(v)}}.
\end{eqnarray*}
Now if
\beq{h}
\prod_{v\in P_u} \frac{1}{1-f(v)/\deg^+(v)}\le 1+h\qquad \forall u \in
X_\ell,
\eeq
then $\sum_{u\in X_\ell}\a_{\ell,u}=1-o(1)$ provided $h=o(1)$.
We next prove  we can choose $h=O(1/\sqrt{\log n})$, which
determines our value of $\e_X$.

Similar to the proof of \eqref{just1} of Lemma \ref{treelemma}
we have, with $\om=\sqrt{\log n}$ that
\beq{just}
\Pr\brac{ \sum_{v \in xPu} f(v) \ge \frac{np}{\om}} \le
\Pr\brac{Bin(n^{2/3+o(1)},p)\ge \frac{np}{\om}}
=O(n^{-10}).
\eeq
Using \eqref{just1}, and \eqref{just} it follows that
$$\sum_{v\in P_u}\frac{f(v)}{\deg^+(v)-f(v)}\le \frac{1}{c\om-1}.$$
For $0<x<1$, $(1-x)^{-1}\leq e^{x/(1-x)}$, and so
\begin{eqnarray}
\prod_{v\in P_u}\frac{1}{1-f(v)/\deg^+d(v)} &\le&
\exp\brac{\sum_{v\in P_u}\bfrac{f(v)}{\deg^+(v)-f(v)}}\nonumber\\
&\leq & \exp\bfrac{1}{c\om-1}=1+O\bfrac{1}{\om}.\label{xx}
\end{eqnarray}
There are at most $n$ trees and $n$ paths per tree and so
\eqref{h}, with $\e_X=h=O(1/\om)=O(1/\sqrt{\log n})$, follows from \eqref{xx}.
This completes the proof of \eqref{LSa}.
\proofend

%%%%%%%%%%%%%%%%%%%%%%%%

The next step is to obtain an estimate of
$\sum_{v\in Y_\ell}\b_{\ell,v}$. The proof is inductive,
moving down the tree $T_y$ level by level.
For brevity we write $d^+(u)=\deg^+(u),\; d_T^-(u)=\deg_T^-(u)$ etc.

Let the random variable $W(y, i)$ be defined by

\[
W(y,i)=\sum_{u \in N^-(y)} \sum_{v \in Y_i} \prod_{z \in vPu} \frac{1}{d^+(z)},
\]
where for $v \in Y_i$ the notation means that the  the unique path $vPuy$
from $v$ to $y$ in $T_y$ passes through $u$, and that $vPu$ is written as $v=z_i,...,z_j,...,z_1=u$
in the product term.

Note that
$$W(y,\ell)=\sum_{v\in Y_\ell}\b_{\ell,v}.$$

Define $W^*(y,i)$ by
\[
W^*(y,i)=\sum_{u \in N^-(y)} \sum_{v \in Y_i} d_T^-(v) \prod_{z \in vPu} \frac{1}{d^+(z)},
\]
where for $v \in Y_{\ell}$ we define $d_T^-(v)=1$ so that $W(y,\ell)=W^*(y,\ell)$.
Note that
\[
W^*(y,1)=\sum_{u \in N^-(y)} \frac{d_T^-(u)}{d^+(u)}.
\]

%%%%%%%%%%%%%%%%%%%%%%%%%

We prove the following lemma for a more general value of $\ell$, as it
is also used in our proof of the upper bound.

\begin{lemma}\label{yLem}
Let
\beq{A2y}
\cA_2(y) = \left\{
D_{n,p}: \sum_{v\in Y_\ell}\b_{\ell,v}= (1+O(1/\sqrt{\log n}))
\frac{1}{np}\sum_{u \in N^-(y)} \frac{d_T^-(u)}{d^+(u)}
\right\}.
\eeq
Let $\ell= \eta \log_{np}n$ where $0 < \eta \le 2/3$.
Then under Assumption 1,
$$\Pr(\exists{y \in V} \text{ such that } \neg \cA_2(y))=O(n^{-\g}).$$

\end{lemma}
\proofstart

The lemma is proved inductively assuming $\cE^-_y$ and $ \cE^+_{Y\sm\{y\}}$.
%and noting that  by \eqref{eqI},  $\Pr(\neg \cE^+_{Y\sm\{y\}})=n^{-1/3-\g}$.
We prove the induction for
$2 \le i \le \ell$, where
by assumption $(np)^{\ell}=O(n^{0.67})$.

Let $\E_{[\ul d^+(i)]} W(y , i)$ be the expectation of $W(y,i)$ over $(d^+(v), v \in Y_{i})$,
conditional on all other degrees  $d^+(u)>0, d_T^-(u), u \in Y_{\le i-1}$ being fixed such that
$|Y_{\leq i-1}|\sim d^-(y)(np)^{i-1}\leq n^{0.67}$ which is true \qs\ from Lemma \ref{treelemma}.

For $v \in Y_i$, $d^+(v)$ is distributed as $D(v)\sim 1+ Bin(\s'(v),p)$,
for some $\s'(v) \in I_0=[n-O(n^{0.67}),n]$. Given the values $\s'(v)$ for $ v \in Y_i$, the $D(v)$
are independent random variables.

For $v \in Y_i$, let $vPu$ be written $vwPu$, where $(v,w) \in T_y$.  Then
\[
\E_{[\ul d^+(i)]} \brac{\prod_{z \in vPu} \frac{1}{d^+(z)}}= \E \bfrac{1}{d^+(v)}
\brac{\prod_{z \in wPu} \frac{1}{d^+(z)}},
\]
where given $\cE^+_{Y\sm\{y\}}$, and $\d=\max(n^{-0.33},n^{-\g})$,

\[
\E \bfrac{1}{D(v)}= (1+O(\d)) \frac{1}{np}.
\]
This follows from the identity
\[
\sum_{j=0}^N \frac{1}{j+1} \binom{N}{j} p^j q^{N-j}x^{j+1}= \frac{1}{(N+1)p} (q+px)^{N+1},
\]
obtained by integrating $(q+px)^N$; and from $\Pr(\neg \cE^+_v)=O(n^{-1-\g})$.
Thus
\begin{eqnarray*}
\E_{[\ul d^+(i)]}W(y,i)&=&
(1+O(\d)) \frac{1}{np}
\sum_{w \in Y_{i-1}} \sum_{v \in N_T^-(w)}\brac{\prod_{z \in wPu} \frac{1}{d^+(z)}}\\
&=&
(1+O(\d)) \frac{1}{np}
\sum_{w \in Y_{i-1}} d_T^-(w)\; \brac{\prod_{z \in wPu}\frac{1}{d^+(z)}}\\
&=&
(1+O(\d)) \frac{1}{np} W^*(y,i-1).
\end{eqnarray*}
To obtain a concentration result, let $U(i)=W(y,i)\cdot((1-o(1))c_0np)^{i}$,
we can write $U(i)=\sum_{v \in Y_i} U_v$, where $U_v$ are independent random variables.
Assuming $\cE^+_{Y\sm\{y\}}$ and that Lemma \ref{treelemma}(ii) holds we have $(c(1-o(1))/C_0)^{i} \le U_v \le 1$.

Let $\e_i= \sqrt{3K \log n /(\E U)}$ for some large constant $K$. Then
\[
\Pr(|U(i)-\E U| \ge \e \E U) \le 2 e^{-\frac{\e^2}{3} \E U} = O(n^{-K}),
\]
and so
\[
\Pr(|W(y,i)-\E W| \ge \e \E W)  = O(n^{-K}).
\]
Note that $\E U \ge |Y_i| (c(1-o(1))/C_0)^{i} \ge (c/2)(c np/2C_0)^{i}$.
Thus  $\e_i\leq1/\sqrt{(A \log n)^{i-1}}$ for some $A>0$ constant.
For $i \ge 2$, $\e_{i}=O(1/\sqrt{\log n})$, and thus $\e_{i}=o(1)$.

In summary, with probability $1-O(n^{-K})$,
\[
W(y,i)= (1+O(\d)+O(\e_i)) \frac{1}{np} W^*(y,i-1).
\]

Continuing in this vein, let
 $\E_{[\ul d_T^-(i-1)]} W^*(y , i-1)$ be the expectation
 of $W^*(y,i-1)$ over $(d_T^-(v), v \in Y_{i-1})$,
conditional on all other degrees  $(d^+(u), d_T^-(u), u \in Y_{\le i-2})$ being fixed.
For $v \in Y_{i-1}$, $d_T^-(v)$ is distributed as $B(v)
\sim Bin(\s(v),p)$ conditional on $\cE^-_{Y\sm Y_\ell}$.
Let $1_\cX$ denote the indicator for an event $\cX$, then
\[
\E B(v)= \E(B(v)\cdot 1_{\cE^-_{Y\sm Y_\ell}})+\E(B(v) \cdot 1_{\neg\cE^-_{Y\sm Y_\ell}}),
\]
and, splitting the second event on $\cD$ gives
\[
\E(B(v) \cdot 1_{\neg\cE^-_{Y\sm Y_\ell}})= O(\D_0 n^{-\g})+ O(n n^{-10}).
\]
Thus, given $\cE^-_{Y\sm Y_\ell}$ we have $\E d_T^-(v)=(1+O(\d))np$.

Thus
\[
\E_{[\ul d_T^-(i-1)]} \brac{ d_T^-(v)\prod_{z \in vPu} \frac{1}{d^+(z)}}= (\E d_T^-(v))
\brac{\prod_{z \in vPu} \frac{1}{d^+(z)}},
\]
and
\[
\E_{[\ul d_T^-(i-1)]} W^*(y,i-1)= (1+O(\d)) np \; W(y,i-1).
\]
Using Lemma \ref{treelemma} (ii) and
 arguments similar to above, for $i \ge 3$ with probability $1-O(n^{-K})$
\[
 W^*(y,i-1)= (1+O(\d)+O(\e_{i-1})) np \; W(y,i-1)
 \]
 completing the induction for $i \ge 3$.

 The final step is to use
 \[
 W(y,2) = (1+O(\d)+O(\e_{2})) \frac{1}{np} \; W^*(y,1),
 \]
 and thus \whp
 \begin{eqnarray*}
 W(y,\ell)&=& \prod_{i=2}^{\ell}(1+O(\d)+O(\e_{i}))^2 \frac{1}{np} \; W^*(y,1)\\
 &=&
 \brac{1+O\bfrac{1}{\sqrt{\log n}} }
  \frac{1}{np} \sum_{u \in N^-(y)} \frac{d_T^-(u)}{d^+(u)}.
 \end{eqnarray*}

Thus from \eqref{eqI}
\[
\Pr(\exists{y \in V} \text{ such that } \neg \cA_2(y))=
O(\Pr(\exists{v \in V}: \deg^\pm(v) \not\in I))=O(n^{-\g}).
\]
 \proofend

\begin{corollary} \label{corol-beta}
Provided Assumption 1 holds, let
\beq{cA2}
\cA_2(y) = \left\{ \sum_{v\in Y_\ell}\b_{\ell,v}= (1+O(1/\sqrt{\log n})) \frac{\deg^-(y)}{np}
\right\} ,
\eeq
then
\begin{equation}
\Pr\brac{\exists{y \in V} \text{ such that } \neg\cA_2(y)}=O(n^{-\g}).\label{LSb}
\end{equation}
\end{corollary}
\proofstart
Referring to \eqref{A2y}, under Assumption 1 and $\neg\cD$, then $d_T^-(u)=\deg^-(u) (1-o(1))$
simultaneously for all
$u\in N^-(y)$
with
probability $1-O(n^{-1-\g})$.
Let $\z=1/\log\log\log n$. A vertex  is {\em normal} if at most
$\z_0=\rdup{4/(\z^3d)} $ of its in-neighbours have
out-degrees which are not in the range $[(1-\z)np,(1+\z)np]$,
and similarly for in-degrees.
Let $\cN(y)$ be the event $y$  is  normal.
We observe that
\[
\Pr(\neg\cN(y)\mid \cE^-_y)\leq 2 \sum_{s=c_0np}^{C_0np}\binom{s}{\z_0}(2e^{-\z^2np/3})^{\z_0}=
O(n^{-\Omega(\log\log\log n)}),
\]
where $\cE^-_y$ is given by \eqref{EsetS}, and
thus (see \eqref{eqI})
\begin{equation}\label{normal}
\Pr(\neg (\cN(y){\cap} \cE^-_y) )=O(n^{-1-\g}).
\end{equation}
Now if $y$ is normal then
$$\deg^-(y)\frac{1-\z}{1+\z}-O(\z_0)\leq \sum_{u \in N^-(y)} \frac{\deg^-(u)}{\deg^+(u)}\leq
\deg^-(y)\frac{1+\z}{1-\z}+O(\z_0).$$
\proofend

Recall the definition of $Z(x,y)$,
\begin{equation}\label{zval}
\Zl(x,y)=\sum_{\substack{u\in X_\ell\\v\in
    Y_\ell}}\a_{\ell,u}\b_{\ell,v}\frac{1_{uv}}{\deg^+(u)},
\end{equation}
where $1_{uv}$ is the indicator for the existence of the edge $(u,v)$
and we take $\frac{1_{uv}}{\deg^+(u)}=0$ if $\deg^+(u)=0$.
The next lemma  gives a high probability  bound for $\Zl(x,y)$.

\begin{lemma}\label{xlem2}
Let
\[
\cA_3(x,y)=\set{Z(x,y)=(1+O(\e_Z)) \frac{\deg^-(y)}{m}},
\]
where $\e_Z=1/(\sqrt{\log n})$. Then given Assumption 1,
\beq{qs3}
\Pr(\exists{x,y} : \neg\cA_3(x,y))=O(n^{-\g}).
\eeq
\end{lemma}
\proofstart

Let
$$\cB=\cB(x,y)=(\cE^+_{X\setminus X_\ell}\cap \cE^-_{Y\setminus Y_\ell}
\cap \cA_1(x,y)\cap \cA_2(y)\cap \cL),$$
where $\cE$  is given by \eqref{EsetS},  $\cA_1,\cA_2$  by
\eqref{cA1}, \eqref{cA2}, and $\cL$ is the event
that Lemma \ref{treelemma} holds.

Let $u \in X_\ell$ and
let $w \in Y \sm Y_\ell$.
As $X_{\ell}\cap Y=\es$,
we know that $u$  is not an in-neighbour
of $w$.
Other
out-edges  of $u$ are unconditioned by the construction of
$\Tl_x,\Tl_y$. Given $Y \sm Y_{\ell} \le n^{2/3+o(1)}$, the
distribution of $\deg^+(u)$ is $Bin(\n,p)$ for some $n-n^{0.67}\leq \n\leq n-1$. Thus
\begin{equation}\label{E1C}
\E\brac{\frac{1_{uv}}{\deg^+(u)}\bigg|\ \cB}=
\sum_{k=1}^{\n}\binom{\n}{k}p^k(1-p)^{\n-k}\frac{k}{\n}\frac{1}{k}=
\frac{1}{n}\brac{1+O(n^{-{0.33}})}.
\end{equation}
Here $k/\n$ is the conditional probability that edge $(u,v)$ is present,
given that $u$ has $k$ out-neighbours.

We use the notation  $\Pr_{\cC}(\cdot)=\Pr(\cdot\mid\cC)$ etc, for any event $\cC$.
From \eqref{cA1}, \eqref{cA2}, \eqref{E1C},
\begin{equation}\label{EZ}
\E_\cB(\Zl)= (1+O(\e_Z))\frac{\deg^-(y)}{m}.
\end{equation}

Conditional on $\cB$, $|Y_{\ell}|\leq n^{2/3+o(1)}$ by construction, and as the edges from
$u$ to $Y_\ell$ are unexposed,
\beq{20}
\Pr_\cB \brac{|N^+(u) \cap Y_{\ell}| \geq 1000} \le
\Pr(Bin(n^{2/3+o(1)},p)\geq 1000)
\leq n^{-10}.
\eeq
Let
$$
\cF=\cF(x,y)=\set{|N^+(u) \cap Y_{\ell}| < 1000,\,\forall u\in X_\ell},
$$
and let $\cG(x,y)=\cB(x,y)\cap \cF(x,y)\cap \cE^+_{X_\ell}$.
The quantity of interest to us is the value of $Z(x,y)$ conditional on $\cG(x,y)$.
We first obtain $\E_\cG(Z)$ from $\E_\cB(Z)$
using
\beq{f1}
\E_{\cB}(Z)= \E_{\cB}( Z\cdot 1_{\cF(x,y)\cap \cE^+_{X_\ell}})+
\E_{\cB}( Z \cdot 1_{\neg[\cF(x,y)\cap \cE^+_{X_\ell}]}).
\eeq
The event $\neg[\cF(x,y)\cap \cE^+_{X_\ell}]
\seq \left[\cF(x,y)\cap \neg \cE^+_{X_\ell}\right] \cup   [\neg \cF(x,y)]$.
Using \eqref{zval}, we obtain
\begin{align}
&\E_{\cB}( Z \cdot 1_{\neg[\cF(x,y)\cap \cE^+_{X_\ell}]})\nonumber \\
&=O(\E_\cB(Z) n^{-\g})+O\brac{\frac{1000}{(c_0np)^{2 \ell}}} |X_\ell|\brac{  O(n^{-(1+\g)})+O(n^{-10})}
+O(|X_{\ell}||Y_\ell| )O(n^{-10})\label{f11}\\
&= \E_\cB(Z)\; O( n^{-\g}).\nonumber
\end{align}
To see this, partition the vertices of $X_\ell$ into sets $R,S$, where vertices in $R$
have out-degree in $[c_0np,C_0np]$, and vertices of $S$ do not.
The first term in \eqref{f11} is the contribution to the
first term in the RHS of \eqref{f1} from the vertices in $R$, multiplied by the probability
of $\neg \cE_{X_\ell}$. Assuming $\cF(x,y)$ holds, the second
term in the RHS of \eqref{f11} is the contribution
to the first term in the RHS of \eqref{f1} from the vertices
in $S$. The last term in the RHS of \eqref{f11} is the contribution
is the contribution
to the first term in the RHS of \eqref{f1} in the case where $\neg\cF(x,y)$ holds.

Thus
\[
\E_\cG(Z)=\frac{\E(Z\cdot 1_{\cB}\cdot 1_{\cF(x,y)\cap \cE^+_{X_\ell}})}{\Pr(\cG)}=
\frac{\E_\cB(Z.1_{\cF(x,y)\cap \cE^+_{X_\ell}})\Pr(\cB)}{\Pr(\cG)},
\]
and so
\beq{n13}
\E_\cG(Z)=\E_\cB(Z)(1+O(n^{-\g}))
=(1+O(\e_Z))
\frac{\deg^-(y)}{np} \frac{1}{n}.
\eeq

\ignore{
Thus
\begin{multline}\label{n13}
\E_\cG(Z)=\frac{\E(Z\cdot 1_{cB}\cdot 1_{\cF(x,y)\cap \cE^+_{X_\ell}})}{\Pr(\cG)}=
\frac{\E_\cB(Z.1_{\cF(x,y)\cap \cE^+_{X_\ell}})\Pr(\cB)}{\Pr(\cG)}\\
\E_\cG(Z)=\E_\cB(Z)(1+O(n^{-\g}))
=(1+O(\e_Z))
\frac{\deg^-(y)}{np} \frac{1}{n}.
\end{multline}
}

We now examine the concentration of $(\Zl \mid \cG)$.
Let $A=1000/((1-o(1))c_0np)^{2\ell+1}$.
It follows from Lemma \ref{treelemma}(ii) that
given $\cG$ we have $Z_u\le A$.
Let $\hZl_u=Z_u/A$, then for $u \in X_{\ell}$, the $\hZl_u$ are independent
random variables, and $0\le \hZl_u\le 1$.
Let $\hZl=\sum_{u\in X_\ell} \hZl_u$ and let $\wh \mu = \E_\cG (\hZl)$.
Thus
\beq{n13a}
\wh \mu =  n^{1/3+o(1)}.
\eeq
It follows from \eqref{Cher1} that if $0\leq \th\leq 1$,
$$\Pr_\cG(|\hZl-\wh \mu| \geq \th\wh \mu )\leq 2e^{-\th^2\wh \mu/3}.$$
With $\th=4(np/\wh \mu)^{1/2}$
we find that,
\[
\Pr_\cG(|\hZl-\wh \mu| \ge 4(np \wh \mu)^{1/2})=o(n^{-4}),
\]
and hence that
$$\Pr_\cG(|Z-\E_\cG Z| \ge 4A(np \wh \mu)^{1/2})=o(n^{-4}).$$
 Using \eqref{n13a} we have $4A(np \;\wh \mu)^{1/2}=O(n^{-7/6+o(1)})$,
 and so
$$
\Pr_\cG\brac{|Z-\E_\cG(Z)|=O\bfrac{1}{n^{7/6+o(1)}}}=1-o(n^{-4}).
$$
We see from \eqref{n13} that
$\E_\cG(Z)= (1+O(\e_Z))\frac{\deg^-(y)}{m}$.
Thus
\beq{gud1}
\Pr_\cG\brac{Z(x,y)  \neq (1+O(\e_Z))\frac{\deg^-(y)}{m}}=o(n^{-4}).
\eeq
%  using \eqref{tilde}
%\beq{gud1}
%\Pr\brac{\brac{\Zl(x,y) \neq (1+O(\e_Z))\frac{\deg^-(y)}{m}}\cap\cG(x,y)}=o(n^{-4}).
%\eeq
Using \eqref{cE}, \eqref{LSa}, \eqref{LSb} and \eqref{20},
\begin{align}
&\Pr\brac{\bigcup_{x,y}\neg \cG(x,y)}\nonumber\\
&\le \Pr(\neg\cE_V)+\Pr(\neg \cL)+\Pr\brac{\bigcup_{x,y}\neg \cF(x,y)}+
\Pr\brac{\bigcup_{x,y}\neg \cA_1(x,y)}+\Pr\brac{\bigcup_y \neg\cA_2(y)} \nonumber\\
&=O(n^{-\g}).\label{notG}
\end{align}
Thus finally, from \eqref{gud1} and \eqref{notG}
\beq{avalu}
\Pr(\exists{x,y} : \neg\cA_3(x,y))= O(n^{-\g}).
\eeq
\proofend

\subsection{Properties needed for an upper bound on the stationary distribution}\label{ubd}
We remind the reader that $np\leq n^\d$.
We first show that small sets of vertices are sparse \whp.
\begin{lemma}\label{smallsets}
Let $\z$ be an arbitrary positive constant. For all
$S\subseteq V,\;|S|\leq s_0=(1-2\z)\Lambda$, $S$ contains at most $|S|$ edges \whp.
\end{lemma}
\proofstart
The expected number of sets $S$ with more than $|S|$ edges can be
bounded by
\begin{eqnarray*}
\sum_{s=3}^{s_0}\binom{n}{s}\binom{s^2}{s+1}p^{s+1}&\leq& \sum_{s=3}^{s_0}(e^2np)^ssep\\
&\le& \exp \brac{-\z\ooi \log n}=o(n^{-\z/2}).
\end{eqnarray*}
\proofend

Let
$$\Lambda=\log_{np}n.$$
We will use the following values:
$$
\begin{array}{lll}
\ell_0=(1+\eta)\Lambda, &\ell_1=(1-10\eta)\Lambda, &\ell_2=11\eta\Lambda.
 %, &\ell_4=\eta\Lambda/20, &\ell_5=9\eta\Lambda/10.
\end{array}
$$
For the upper bound we need to slightly alter our definition of breadth-first trees
and call them $\Tu_x,\Tu_y$.
This time we grow $\Tu_x$ to a depth $\ell_1$ and $\Tu_y$
to a relatively small depth $\ell_2$.
With this choice, Lemma \ref{smallsets} implies that $Y$ will
contain no more than $|Y|$ edges \whp.
This reduces the complexity of the argument.
We fix $x,y$ and grow  $\Tu_x$ from $x$ to a depth $\ell_1$,
and $\Tu_y$ into $y$ to a depth $\ell_2$.
The definition of $\Tu_x$ is slightly different from $\Tl$, but we retain some of the notation.

{\bf Construction of $\Tu_x$.}
We build a tree $\Tu_x$, much as in Section \ref{lbd},
by growing a breadth-first out-tree from $x$ to depth $\ell$.
The difference is that we construct $\Tu_x$ before $\Tu_y$, so that
$\Tu_x$ is not disjoint from $Y$.
As before, let $X_0=\{x\}$, and $X_i,\; i \ge 1$ be the $i$-th level set of the tree.
Let $\Tu_x(i)$ denote the BFS tree up to and including level $i$, and let $\Tu_x=\Tu_x(\ell_1)$.
Let $X_{\leq i} = \cup_{j \le i} X_j$, and let ${\sX}=X_{\leq \ell_1}$.
In Section \ref{upperbound} below we
will need to consider a larger set $X_{\leq\ell_3}$ where $\ell_3=(1-\eta/10)\La$.

{\bf Construction of $\Tu_y$.}
Our upper bound construction  of $\Tu_y$ is
the same
as for the lower bound, except that we only grow it to depth $\ell_2$.

Our aim is to prove an upper bound similar to the lower bound proved in Lemma \ref{xlem2}.
For $u\in X_{\ell_1}$ we let
$$\a_{\ell_1,u}=\Pr(\cW_x(\ell_1)=u)$$
where
\begin{equation}\label{upperalpha}
\sum_{u\in {\sX}}\a_{\ell_1,u}\leq 1.
\end{equation}
The RHS of \eqref{upperalpha} is one, except when we fail to construct $\Tu_y$ to level $\ell_2$.

This is the only place where we write a structural property of $D_{n,p}$ in terms of
a walk probability. This is of course valid, since $\a_{\ell_1,u}$ is the sum over walks
of length $\ell_1$ from $x$ to $u$ of the product of reciprocals of out-degrees.
Fortunately, all we need is \eqref{upperalpha}.

We also define the $\b_{i,v}$ as we did in \eqref{betaiv} and now we let
\begin{equation}\label{Giveny}
\Zu(x,y)=Z^{up}(x,y)=\sum_{ \substack{u\in {\sX}\\v\in
    Y_{\ell_2}\setminus {\sX}}}
    \a_{\ell_1,u}\b_{\ell_2,v}\frac{1_{uv}}{\deg^+(u)}.
\end{equation}
The following lemma follows from Corollary \ref{corol-beta}.
\begin{lemma}\label{xlem4}
Let $\ell$ be as in \eqref{ell}.
If $2\leq k\leq \ell$ then for some $\e_Y=o(1)$ we have
$$\Pr\brac{\sum_{v\in Y_{k}}\b_{k,v}\geq (1+\e_Y)\frac{\deg^-(y)}{np}}=o(n^{-1-\g/2})$$
where $\g$ is as in \eqref{eqI}.
\end{lemma}

It follows by an argument similar to that for Lemma \ref{xlem2} that
\begin{lemma}\label{xlem5}
\begin{equation}\label{Zconcx}
\Pr\brac{\exists x,y:\;\Zu(x,y)\geq (1+o(1))\frac{\deg^-(y)}{m}}=O(n^{-\g/2}).
\end{equation}
\end{lemma}
In computing the expectation of $\Zu$, some of the vertices in   ${\sX}$ of $\Tu_x$ may be
inspected in our construction of $\Tu_y$, or of $\Tu_x$ up to level $\ell_1$.
Thus $\E(1_{uv}/\deg^+(u))\leq (1/n) \ooi$, (see \eqref{E1C}).

\begin{Remark}\label{rem2}
The upper bound for $\Zu(x,y)$
obtained above is parameterized by $\ell_0=(1+\eta)\La$.
Provided $\eta >0$ constant, so that Lemma \ref{xlem5} holds, we can apply this argument
simultaneously for $n^{\g/3}$ different values of $\eta$.
\end{Remark}

We next prove a lemma about non-tree edges inside ${\sX}$,
and  edges from ${\sX}$ to $Y\sm Y_{\ell_2}$.
\begin{lemma}\label{back}\
\begin{description}
\item[(a)]
Let $\ell_3=(1-\eta/10)\Lambda$ and
\[
\cLa(\ell_3)= \set{\forall \; z\in X_{\leq \ell_3}:\;z\ \text{ has } \le
100/\eta \text{ in-neighbours  in } X_{\leq \ell_3}}.
\]
Then $\Pr(\neg\cL_a(\ell_3))=O(n^{-9}).$
\item[(b)] Let
$X^\circ_\ell=\set{v \in X_\ell: N^+(v) \cap X_{\leq\ell} \ne \es}$ and
\[
\cL_b(\ell)=\set{|X^\circ_\ell| \le 18{\D_0}^{2\ell}p+\log^2 n}.
\]
Then $\Pr(\neg\cL_b(\ell))=O(n^{-10})$ for $\ell\leq\ell_3$.
\item[(c)] Let
\beq{t0}
t_0=\rdup{\frac{K\La}{\log np}}\qquad\text{where $K=2\log(100C_0/\eta c_0)$}.
\eeq
Fix $t\leq t_0$ and $i>2\eta\La$ and let
$$S_{i,t}^\circ=\set{z\in {\sX}:\;z\text{ is reachable from $X_i^\circ$ in at most $t$ steps}}.$$
Then let $A^\circ=A^\circ(x,y,t)$ be the number of edges from
$S^\circ_{i,t}\cap {\sX}$ to $Y_{\ell_2-t}$. Then,
$$\Pr(\exists x,y,t:\;A^\circ\geq \log n)=O(n^{-10}).$$
\item[(d)]
Let $A=A(x,y,t)$ be the number of edges between $X_{\ell_1}$ and $Y_{\ell_2-t}\sm {\sX}$,
where $t_0<t\leq \ell_2-1$.
\[
\Pr(\exists x,y,t:\;A \ge 9 |X_{\ell_1}||Y_{\ell_2-t}|p + \log^2 n) =o(n^{-10}).
\]
\end{description}
\end{lemma}
\proofstart\\
(a)
Let $z \in X(\ell_3)$.
Let $\z$ be the number of
in-neighbours of $z$ in $X_{\leq \ell_3}$.
In the construction of $\Tu_x(\ell_3)$, we only exposed one in-neighbour of $z$.
Thus $\z$ is distributed as $1+Bin(|X_{\leq\ell_3}|,p)\leq 1+Bin(\D_0^{\ell_3},p)
+n\Pr(\cD)$.
We apply \eqref{Cher2} and  \eqref{largedegree} to deal
with the binomial.
Hence if $r+1=100/\eta$,
$$\Pr(\z\geq r+1)\leq \D_0^{r\ell_3}p^r+n^{-10}e^{-10np}\leq 2n^{r(\d-\eta/10)}+n^{-10}e^{-10np}=
O(n^{-9}).$$
Part (a) of the
lemma follows.

(b) For $v \in X_\ell$ the out edges of $v$ are unconditioned during the construction
of $\Tu_x(\ell)$.  The number of out edges of
$v$ to $X_{\leq\ell}$ is  $Bin(|X_{\leq\ell}|,p)$.
Unless $\cD$ occurs, $|X_{\leq\ell}| \le 2{\D_0}^{\ell}$ and
\[
\Pr(|N^+(v) \cap X_{\leq\ell}|> 0\; \mid  \neg \cD) %|X_{\leq\ell}| \le 2{\D_0}^{\ell})
\le 1-(1-p)^{2{\D_0}^{\ell}}\le2{\D_0}^{\ell}p,
\]
and
$$\E( |X^\circ_\ell|\; \mid \neg \cD) %X_{\leq\ell}| \le 2{\D_0}^{\ell})
 \le 2{\D_0}^{2\ell}p.$$
By \eqref{Cher2}
\[
\Pr(|X^\circ_\ell| \ge 18{\D_0}^{2\ell}p+\log^2 n)=O(n^{-10})+\Pr(\cD)=O(n^{-10}).
\]
(c)
Let $S(u,t')$ be the set of vertices in ${\sX}$ that a walk starting from
$u\in X_i$ can reach in $\ell_1-i+ t-t'$ steps. Thus unless $\cD$ occurs,
$|S(u,t')| \le {\D_0}^{\ell_1-i+ t-t'}$.
So, given $\neg\cD$,
\beq{Scirc}
|S^{\circ}_{i,t}| \le 2|X^{\circ}_i| {\D_0}^{\ell_1-i+ t}.
\eeq
We can assume that, after constructing $\Tu_x$ we construct
$\Tu_y$ to level $Y_{\ell_2-t}$, and then inspect the edges
from $S^{\circ}_{i,t}$ to $Y_{\ell_2-t}\sm {\sX}$. These edges are
unconditioned at this point and their number $A$ is stochastically dominated by
$Bin(|S^{\circ}_{i,t}|\;|Y_{\ell_2-t}|, p)$ .
Given $\cL_b(i)$ of part (b) of this lemma,
\beq{Xcirc}
|X^{\circ}_i| \le 18{\D_0}^{2i}p+\log^2 n.
\eeq
Let $i=a \La$, where $2\eta \le a \le 1-10\eta$.

{\bf Case  $2\eta \le a \le (1+\e)/2$ for some small $\e>0$ constant.}\\
Using \eqref{Scirc}, \eqref{Xcirc} and $|Y_{\ell_2-t}| \le {\D_0}^{\ell_2-t}$ gives
\begin{eqnarray*}
\E A^\circ &\le & (18{\D_0}^{2i}p+\log^2 n) 2 {\D_0}^{\ell_1-i+t} {\D_0}^{\ell_2-t} p
+n^2(\Pr(\cD)+\Pr(\cL_b(i)))\\
&\le& 36 C_0^{\ell_0+\ell_1}\; p^2 \;(np)^{\ell_0+i} +
2 C_0^{\ell_0}(\log^2n)\; p \;(np)^{\ell_0-i}+O(n^{-9})\\
&\le& 36 C_0^{\ell_0+\ell_1} (np)^2 n^{-\half+\eta+\e} +2C_0^{\ell_0}\log^2n (np) n^{-\eta}
+O(n^{-9})\\
&=&O(n^{-\eta/2}).
\end{eqnarray*}
{\bf Case  $(1+\e)/2\le a \le 1-10\eta$.}\\
For $i \ge (1+\e)/2 \La$, $|X^{\circ}_i| \le 20{\D_0}^{2i}p$. Thus
\begin{eqnarray*}
\E A^\circ &\le & 20{\D_0}^{2i}p \cdot  2 {\D_0}^{\ell_1-i+t} {\D_0}^{\ell_2-t} p
+n^2(\Pr(\cD)+\Pr(\cL_b(i)))\\
&\le& 40 C_0^{\ell_0+\ell_1}\; p^2\; (np)^{\ell_0+i}+O(n^{-9})\\
&\le& 40 C_0^{\ell_0+\ell_1} (np)^2 n^{-9\eta}+O(n^{-9})\\
&=&O(n^{-\eta/2}).
\end{eqnarray*}

In either case, with probability $1-o(n^{-10})$, $A \le \log n$.

(d) After growing $\Tu_x$ to level $\ell_1$, we grow $\Tu_y$ to level $\ell_2-t$.
Then $A(t)$ has a binomial distribution and $\E A(t) \le |X_{\ell_1}||Y_{\ell_2-t}|p$.
The result follows from the Chernoff inequality.
\proofend

\subsection{Small average degree: $1+o(1)\leq d\leq 2$}
This section contains further lemmas needed for the case $1+o(1)\leq d\leq 2.$

We will assume now that
$$1+o(1)\leq d\leq 2.$$
Let a vertex be
{\em small} if it has in-degree or out-degree at most $np/20$ and
{\em large} otherwise. Let {\em weak distance} refer to
distance in the  underlying {undirected} graph of
$D_{n,p}$.
\blem{newlem0}\
\begin{description}
\item[(a)] {\bf Whp} there are fewer than $n^{1/5}$ small vertices.
\item[(b)] If $np\geq 2\log n$ then \whp\ there are no small vertices.
\item[(c)] {\bf Whp} every pair of small vertices are at weak distance at least
$$\ell_{10}=\frac{\log n}{10\log\log n}$$
apart.
\item[(d)] {\bf Whp} there does not exist a vertex $v$ with
  $\max\set{\deg^+(v),\deg^-(v)}\leq \log n/20$.
\item[(e)]
Let $\varsigma^*(v)$ be given by \eqref{varsig}.
{\bf Whp} for all vertices $y$,
$$\sum_{u\in
    N^-(y)}\frac{\deg^-(u)}{\deg^+(u)}=(1+o(1))(\deg^-(y)+\varsigma^*(y)).$$
\end{description}
\elem
\proofstart\\
(a) The expected number of small vertices
is at most
\beq{qwerty}
n\sum_{k=0}^{\log n/20}\binom{n-1}{k}p^kq^{n-1-k}=O(n^{.1998}).
\eeq
Part (a) now follows from the Markov inequality.

(b) For $np\geq 2\log n$ the RHS of \eqref{qwerty} is $o(1)$.

(c) The expected number of pairs of small vertices at distance $\ell_{10}$ or
less is at most
\begin{multline*}
n^2\sum_{k=0}^{\ell_{10}}2^kn^kp^{k+1}\brac{2\sum_{l=0}^{\log
    n/20}\binom{n-1}{l}p^lq^{n-1-l}}^2=\\
O(n \ell_{10} (2d\log
n)^{\ell_{10}+1}(20ed)^{\log n/10}n^{-2d})=O(n\cdot n^{1/10+o(1)}\cdot n^{1/2}\cdot n^{-2})=o(1).
\end{multline*}
(d) The expected number of vertices with small out- and in-degree is
$O(n^{1-2\times .8002})=o(1)$.

(e)
For $1\leq k\leq  {\D_0}$  let
$$\l_k=\begin{cases}1&1\leq k\leq \frac{\log n}{(\log\log n)^4}\\(\log\log
  n)^4&\frac{\log n}{(\log\log n)^4}\leq k\leq \D_0 \end{cases}.$$
Let $\e=\frac{1}{\log\log n}$.
The probability {that} there exists a vertex of in-degree $k\in
[1,\D_0]$ with
$\l_k$ in-neighbours of {in or out-}degree outside $(1\pm\e)np$, is
bounded by
$$\sum_{k=1}^{\D_0}n\binom{n-1}{k}p^kq^{n-1-k}\binom{k}{\l_k}
(4e^{-\e^2np/3})^{\l_k}\leq\sum_{k=1}^{\D_0}
  2n^{1-d}\brac{\frac{nep}{k}\cdot 2 \cdot
    n^{-\e^2d\l_k/(4k)}}^k=o(1).$$
Now assume that there are fewer than $\l_k$ neighbours of $v$ of
{in or out-}degree outside $(1\pm\e)np$.
Assuming at most one neighbour $w$ of $y$
is small,
$$\sum_{u\in
    N^-(y)\setminus \set{w}}\frac{\deg^-(u)}{\deg^+(u)}=\begin{cases}(1+O(\e))k&1\leq
    k\leq \frac{\log n}{(\log\log n)^4}\\(1+O(\e))
    (k-\l_k)+O(\l_k)&\frac{\log n}{(\log\log n)^4}\leq k\leq \D_0
\end{cases}.$$

This completes the proof of the lemma.
\proofend

  Let {\em weak} distance refer to distance in the underlying graph of $D_{n,p}$,
and let a cycle in the underlying graph  be called a {\em
  weak} cycle.

\blem{newlem2}
{\bf Whp} there does not exist a small vertex that is within weak distance
$\ell_{10}$ of a weak cycle $C$ of length at most $\ell_{10}$.
\elem
\proofstart
Let $v,C$ be such a pair. Let $|C|=i$ and $j$ be the weak distance of $v$ from $C$.
The probability that  such a pair exists is at most
\begin{multline*}
\sum_{i=3}^{\ell_{10}}(2np)^ii\sum_{j=0}^{\ell_{10}}(2np)^j\sum_{l=0}^{\log
    n/20}2\binom{n-1}{l}p^l q^{n-1-l}\\
=O(n^{1/10+o(1)}\cdot n^{1/10+o(1)}\cdot n^{-4/5+o(1)})=o(1).
\end{multline*}
\proofend

\ignore{
\begin{lemma}\label{lemcov}
Let $V^*=\set{v\in V:\;\deg^-(v)=k^*=\rdup{(d-1)\log n}\ and\ \deg^+(v)=\rdup{d\log n}}$.
Then \whp\ there is no weak cycle of length less than $\La/2$ that contains a vertex
of $V^*$.
\end{lemma}
\proofstart
The probability that there is such a cycle can be bounded by
\begin{align*}
&\sum_{r=3}^{\La/2}\binom{n}{r}rp^r\Pr(Bin(n-2,p)=k^*-1)\\
&\leq \sum_{r=3}^{\La/2}r(np)^r\binom{n-2}{k^*-1}p^{k^*-1}(1-p)^{n-k^*-1}\\
&\leq \La n^{1/2-d}e^{k^*p}\bfrac{nep}{k^*-1}^{k^*-1}\\
&\leq \La n^{-1/2}e^{k^*p}n^{(d-1)\log(d/(d-1))}
\end{align*}
\proofend
}
\section{Analysis of the random walk: Estimating the stationary distribution}\label{SteadyState}
In this section we keep Assumption 1 and assume that we are dealing with a digraph
 which has all of the high probability properties of the previous section.

\subsection{Lower Bound on the stationary distribution}\label{lowerbound}
We use the properties described in Section \ref{lbd}.
We derive a lower bound on $P_x^{2\ell+1}(y)$.
For this lower bound we only consider $(x,y)$-paths of length $2\ell+1$
consisting of a $\Tl_x$ path from $x$ to $X_{\ell}$ followed by an edge
from $X_\ell$ to $Y_\ell$ and then a $\Tl_y$ path to $y$.
The probability of following such a path is $\Zl(x,y)$,
see \eqref{Given}. Lemma \ref{qs3} implies that

\beq{lowpi}
P^{(2\ell+1)}_x(y)\geq (1-o(1))\frac{\deg^-(y)}{m}\mbox{ for all }v\in V.
\eeq

\blem{kl}
For all $y\in V$,
$$\p_y\geq (1-o(1))\frac{\deg^-(y)}{m}.$$
\elem
\proofstart
It follows from \eqref{lowpi} that for any $y\in V$,
\begin{equation}\label{finally}
\p_y=\sum_{x\in V}\p_xP^{(2\ell+1)}_x(y)\geq
(1-o(1))\frac{\deg^-(y)}{m}\sum_{x\in
  V}\p_x=(1-o(1))\frac{\deg^-(y)}{m}.
\end{equation}
\proofend

\subsection{Upper Bound on the stationary distribution}\label{upperbound}

Lemma \ref{kl} above proves that the expression in Theorem \ref{lemsteady} is
a lower bound on the stationary distribution.
As $\sum \p_y=1$, this can be used to derive an
upper bound of $\p_y\leq (1+o(1))\frac{\deg^-(y)}{m}$ which holds for all but $o(n)$ vertices $y$.
In this section we
extend this upper bound {\em to all} $y\in V$.

We use the properties described in Section \ref{ubd}.
We now consider the probability of various types of walks of length $\ell_0+1$ from $x$ to $y$.
Some of these walks are simple directed paths
in  BFS trees {constructed in a similar way to}  the lower bound,
and some use back edges of these BFS trees, or contain cycles etc.
We will upper bound $P_x^{\ell_0+1}(y)$
as a sum
\begin{equation}\label{Pxya}
P_x^{\ell_0+1}(y)\leq Z_x^{\ell_0+1}(y)+S_x^{\ell_0+1}(y)+Q_x^{\ell_0+1}(y)+R_x^{\ell_0+1}(y),
\end{equation}
where the definitions of the probabilities on the right hand side are described below.

\begin{description}
\item[$\ul{Z_x^{\ell_0+1}(y)}$.]
This is the probability that $\cW_x(\ell_0+1)=y$ and the $(\ell_1+1)$th edge $(u,v)$ is such that
 $u \in {\sX}$ and  $v\in Y_{\ell_2}\sm {\sX}$, and
the last $\ell_2$ steps of the walk use edges of the tree $\Tu_y$.
These are the simplest walks to describe.
They go through $\Tu_x$ {for $\ell_1$ steps} and
then level by level through $\Tu_y$.
They make up almost all of the walk probability.

\item[$\ul{S_x^{\ell_0+1}(y)}$.] This is the probability
that $\cW_x(\ell_0+1)=y$ goes from $x$ to $y$
without leaving ${\sX}$.
This includes any special cases such as, for example,  a walk $xyxy...xy$
based on the existence of a cycle $(x,y), (y,x)$
in the digraph.

\ignore{
and (i) does not leave ${\sX}$ until on or
after step $\ell_1+t_0$, ($t_0$ defined in \eqref{t0}) and (ii) on leaving
${\sX}$ follows the unique path in $\Tu_y$ to $y$.}

\item[$\ul{ Q_x^{\ell_0+1}(y)}$.] This is the probability that $\cW_x(\ell_0+1)=y$ and
the $(\ell_1+1)$th edge $(u,v)$ is such that $v\in Y_{\ell_2}\cap  {\sX}$ and
the last $\ell_2$ steps of the walk use edges of the tree $\Tu_y$.
We exclude
walks within ${\sX}$ that are counted in $\ul{S_x^{\ell_0+1}(y)}$.
%So we only count walks that leave ${\sX}$ by step before $\ell_1+t_0$.

\item[$\ul{R_x^{\ell_0+1}(y)}$.] This is the probability that $\cW_x(\ell_0+1)=y$ and during
the last $\ell_2$ steps, the walk uses some edge which is a back or cross edge
with respect to the tree $\Tu_y$.

\end{description}

\underline{\bf Upper bound for $Z_x^{\ell_0+1}(y)$.}

It follows from \eqref{Zconcx} that
\beq{upZ}
Z_x^{\ell_0+1}(y)\leq (1+o(1))\frac{\deg^-{y}}{m}.
\eeq
\underline{\bf Upper bound for $S_x^{\ell_0+1}(y)$.}

Let $\cW_x$ be a walk of length $t$ in ${\sX}$, and let $\cW_x(t)=v$.
Let $d^{-}_{\max}=\max_{w \in {\sX}}|N^-(w) \cap {\sX}|$. Tracing back from $v$ for
$t$ steps, the number of walks length $t$ in ${\sX}$ terminating at $v$ is at most $(d^-_{\max})^t$;
so this serves as an upper bound on the number of walks from $x$ to $v$ of this length.
By Lemma \ref{back}(a), we may assume that  $d^-_{\max} \le 100/\eta$.

Applying this description,
there can be at most $(100/\eta)^{\ell_0+1}$  walks of length $\ell_0+1$
 from $x$ to $y$, which do not exit from ${\sX}$. We conclude that
\beq{Sval}
S_x^{\ell_0+1}(y)\le \bfrac{(100/\eta)}{c_0np}^{\ell_0+1}= o\bfrac{1}{n^{1+\eta/2}}.
\eeq

\underline{\bf Upper bound for $Q_x^{\ell_0+1}(y)$.}

We say that a walk $\cW_x$ {\em delays for $t$ steps},
if $\cW_x$ {\em exits ${\sX}$ for the first time} at step $\ell_1+t$.
A walk {\em delays} at level $i$, if the walk takes
a cross edge (to the same level $i$) or a back edge (to a level $j<i$)
i.e. a {\em non-tree edge} $e=(u,v)$ contained in $X$ that is not part of $\Tu_x$.

\begin{lemma}
Let $t_0=\rdup{\frac{K\La}{\log np} }$ where $K=2\log(100C_0/\eta c_0)$, then
\[
\Pr(\cW_x(\ell_0+1)=y \text{ and } \cW_x \text
{ delays for } t_0 \text{ or more steps})=o(1/n).
\]
\end{lemma}
\proofstart
The only way for a walk to exit from ${\sX}$ is via $X_{\ell_1}$
(recall that edges oriented out from $X_i$ end in $X_{i+1}$).
Let $\cW_x$ be an $(x,y)$-walk which delays for $t$ steps, and then
takes  edge $e=(u,v)$  between $X_{\ell_1}$ and $Y_{\ell_2-t}\sm {\sX}$. There are at most
$(100/\eta)^{\ell_1+t}$ walks  of length $\ell_1+t$
from $x$ to $u$ within ${\sX}$. After reaching vertex $v$, $\cW_x$ follows the
 unique path from $v$ to $y$ in $T_Y$. Applying Lemma \ref{back}(d) we see that
the total probability  $P^{\dag}(t)$ of such $(x,y)$-walks of length $\ell_0+1$ and delay $t$ is
\begin{eqnarray*}
P^{\dag}(t)&\le&\frac{(100/\eta)^{\ell_1+t}
\brac{9 |X_{\ell_1}||Y_{\ell_2-t}|p + \log^2 n}}{(c_0np)^{\ell_0+1}}\\
&\le&
\frac{(100/\eta)^{\ell_1+t}\brac{9C_0^{\ell_0}(np)^{\ell_0-t}p +
\log^2n}}{(c_0np)^{\ell_0+1}}\\
&=&O \brac{\frac{1}{n}\; \bfrac{(100/\eta)C_0}{c_0}^{\ell_0}\brac{\frac{1}{(np)^{t}}+
\frac{\log^2n}{n^{\eta}}}}.
\end{eqnarray*}
So,
\beq{Pdag1}
P^{\dag}(\ge t_0)=\sum_{t \ge t_0} P^{\dag}(t)
=O \brac{\frac{A^{\La}}{n}\brac{\frac{1}{(np)^{t_0}}+\frac{\log^2n}{n^{\eta}}}},
\eeq
where $A=(100C_0/\eta c_0)^{1+\eta}$.

Now $A^{\La}=n^{o(1)}$. Also
$A^{\La}/np=o(1)$ if $\log^2np\geq 2(\log A)(\log n)$ in which case the RHS of \eqref{Pdag1}
is $o(1/n)$, which is what we need to show. So assume now that $\log^2np\leq 2(\log A)(\log n)$.
This means that $\La\to\infty$ and then
$$\frac{A^{\La}}{(np)^{t_0}}\leq \bfrac{A}{e^K}^{\La}\to 0.$$
Thus in both cases
\beq{Pdag}
P^{\dag}(\ge t_0)=o(1/n).
\eeq
\proofend

We can now focus on walks with delay $t$, where $1 \le t < t_0$.
A {\em non-tree} edge of $X$ is an edge induced by $X$ which is
not an edge of $\Tu_x$.
For $4i \le (1-\e) \La$,
 Lemma \ref{smallsets} implies  that \whp\ the set $U=X_{\leq i}$
contains at most $|S|$ edges.
 For, if $U$ contained more
than $|U|+1$ edges then it would contain two distinct cycles
$C_1,C_2$. In which case, $C_1,C_2$ and the shortest undirected path in $U$ joining them
would form a set $S$ which satisfies the
conditions of Lemma \ref{smallsets}.
Thus there is at most one {\em non-tree edge} $e=(u,v)$
 contained in $X_{\leq (1-\e)\La/4}$.

Let $\th = 2\eta \La$.
We classify walks into two types.
\begin{description}
\item[ Type 1 Walks.]
These
have a delay caused by using a non-tree edge of $X_{\leq\th}$, but no delay arising
at any level $i > \th$. Thus, once the walk finally exits
$X_{\th}$ to $X_{\th+1}$ it moves forward at each step
towards $X_{\ell_1}$, and then exits
to $Y_{\ell_2-t}\sm {\sX}$.

\item[Type 2 Walks.] These have a delay arising at some
level $X_i, \; i > \th$. We do not exclude  previous delays
occurring in $X_{\leq\th}$,
or subsequent delays at any level.
\end{description}

\underline{Type 1 Walks}.
We can assume that
$X_{\leq\th}$ induces exactly one non-tree edge $e=(u,v)$.
Let $u \in X_i$ then $ v \in X_j, \; j \le i$.
There are two cases.

(a) {\bf $e$ is a cross edge, or back edge
not inducing a directed cycle}.
%The delay is $t=1+(i-j)$, $a=i$, $w=v$.

Here the delay is $t=i+1-j$ and this is less than $t_0$ by assumption. Then, as we will see,
\beq{T1a}
\Pr(\text{Type 1(a) walk})\leq \frac{1}{(c_0np)^2}\frac{1}{\deg^+(w)}\sum_{w\in N^+(v)}
Z_w^{(\ell_0-i+j+1)}(y)=O\bfrac{1}{n(np)^2}.
\eeq
The term $1/(c_0np)^2$ arises from the walk having to take the out-neighbour of $x$ that leads to
$u$ in $\Tu_x$ and then having to take the edge $(u,v)$. The next step of the walk is to choose
$w\in N^+(v)$ and it must then follow a path to $y$ level by level through the two trees.
The value of $Z_w^{\ell_0-i+j+1}(y)$
can be obtained as follows.
Let $\ell_0'=\ell_0-(i-j)-1$, then as $t < t_0=o(\La)$ we have that $\ell_0'\sim\ell_0$.
For $w \in N^+(v)$
replace $\ell_0$ by $\ell_0'$  in \eqref{Zconcx} above, to obtain
$Z(w,y)=O(1/n)$, see Remark \ref{rem2}. This verifies \eqref{T1a}.

(b) {\bf $e$ is a back edge inducing a directed cycle}.

Let $xPu$ be the path
from $x$ to $u$ in $\Tu_x(\th)$. As $v$ is a vertex of $xPu$,
we can write $xPu=xPv,vPu$ and cycle $C=vCv=vPu,(u,v)$.
Let $\s\ge 2$ be the length of $C$. For some $w$ in $vPu$ the
walk is of the form $xPv,vPw, (wCw)^k,wPz$, where $wCw$
is $C$ started at $w$, the walk goes round $wCw$, $k$ times and exits at $w$ to $u'\in N^+(w)\sm C$
and then moves forward along $wPz$ to $z \in X_{\ell_1}$ and then onto $y$.
The delay is $t=k\s$ and this is less than $t_0$ by assumption

We claim that
\beq{PType1}
\Pr(\text{Type 1(b) walk}) \le \sum_{w\in C}\sum_{k\geq 1}(c_0np)^{-k\s}\frac{1}{\deg^+(w)-1}
\sum_{u'\in N^+(w)\sm C}Z_{u'}^{\ell_0-k\s+1}(y)=O\bfrac{1}{n(np)^2}.
\eeq

The term $(c_0np)^{-k\s}$ accounts for having to go round $C$ $k$ times and we can argue that
$Z_{u'}^{\ell_0-k\s+1}(y)=O(1/n)$ as we did for Type 1(a) walks.

So from \eqref{T1a} and \eqref{PType1} we have that
\beq{Type1}
\Pr(\text{Type 1 walk}) =O\bfrac{1}{n(np)^2}.
\eeq

\underline{Type 2 Walks}.
Suppose $\cW_x$ is a walk which exits ${\sX}$ at
step $\ell_1+t$ and is delayed at some level $i > \th$
by using an edge $(u,v)$.
The walk arrives at vertex  $u \in X_i$ for the first
time at some step $i+t'$ and traverses a cross or back edge to
$v \in X_j, \; j \le i$.

A contributing walk will have to use one of the $A^\circ(x,y,t)\leq \log n$ edges described in
Lemma \ref{back}(c).
By Lemma \ref{back}(a)
there are at most $(100/\eta)^{\ell_1+t}\log n$ from $x$ to $u \in X^{\circ}_i$.
Once the walk reaches $w \in Y_{\ell_2-t}$ there is (by assumption)
a unique  path in $\Tu_y$ from $w$ to $y$.
Let $P(i,t)$ be the probability of these Type 2 walks, then
\beq{Pit}
P(i,t) \le \frac{(100/\eta)^{\ell_1+t} \log n}{(c_0np)^{\ell_0+1}} = O \bfrac{1}{n^{1+\eta/2}}.
\eeq

Thus finally from  \eqref{Pdag}, \eqref{Type1}, \eqref{Pit}
\beq{Qval}
Q_x^{\ell_0+1}(y) = P^{\dag}(\ge t_0)+\Pr(\text{Type 1 walk}) +
\sum_{1 \le t \le t_0}\sum_{\th \le i \le \ell_i} P(i,t)
= O\brac{\frac{1}{n} \frac{1}{(np)^2}}.
\eeq

\underline{\bf Upper bound for $R_x^{\ell_0+1}(y)$.}

Let $Y=Y_{\leq\ell_2}$ be the vertex set of $\Tu_y(\ell_2)$.
We assume that $Y$ induces a unique edge $e=(u,v)$ which is not in $\Tu_y$.
Note that the condition that $|Y|$ induces at most $|Y|$ edges holds,
even if we replace $\ell_2$ with $2 \ell_2$
 based on the construction of $T_{Y}(2\ell_2)$ to depth $2\ell_2$,
by branching backwards from $y$.
We consider two cases.

{\bf  (i) $e$ is a cross or forward edge, or back edge not inducing a directed cycle.}

We
have $u\in Y_i,v\in Y_j$ for some $i\le j \le \ell_2$.
We suppose the $(x,y)$-walk is of the form $xWu,(u,v),vWy$ where $u \not \in vWy$,
so that $vWy$ is a unique path in $\Tu_y$.

\underline{Case  1: $i > (4\eta/5) \La$.}

Let $\ell_3 = (1-\eta/10)\La$.
The length of the path $(u,v),vWy$ is $j$, so the length of the walk $xWu$ is $\ell_0-j+1$.
Let $h$ be the distance from $u$ to $X_{\leq\ell_3}$ in $\Tu_y$.
Then
\[
 h=\max\set{0,\ell_0-\ell_3-j+1} \le \max\set{0,\ell_0-\ell_3-i+1}.
\]

Let $w \in X_{\leq\ell_3}$. By Lemma \ref{back},
the number of $(x,w)$-walks of length $\ell\le \ell_3$ in $X_{\leq\ell_3}$
passing through $w$ at step $\ell$
 is bounded by $(100/\eta)^{\ell_3}$. The
 the number of walks length $h$ from $u$ to
 $ X_{\leq\ell_3}$ is at most ${\D_0}^h$.
Thus, the number of $(x,y)$-walks passing through $e=(u,v)$ is bounded by
$(100/\eta)^{\ell_3}{\D_0}^{h}$.
Thus
\beq{Rval1}
R_x^{\ell_0+1}(y) =
O\bfrac{(100/\eta)^{\ell_3}{\D_0}^{9\eta\La/10}}{(c_0np)^{\ell_0+1}}=O(n^{-1-\eta/20}).
\eeq

\underline{Case 2: $0< i \le (4\eta/5) \La$.}

Let $i=a\La$. Let $\eta'=\eta(1-a)\;, \ell_1'=(1-10\eta')\La,\;\ell_2'=11\eta'\La$, and
let $\ell_0'=\ell_1'+\ell_2'$.
As observed above, the vertex set $U$ of the tree $T_U$ of height $\ell_2'$ above $u$
 induces no extra edges, so we can apply the
upper bound result for walks of length $\ell_0'+1$ from $x$ to $u$ based on the assumption
$R_x^{\ell_0'+1}(u)=0$.
Thus
\[
P_{x}^{\ell_0'+1}(u) \le (1+o(1))
\frac{\deg^-(u)}{m}.
\]
The probability the walk then follows the path $(u,v),vPy$ is $O(1/(np)^2)$.
Thus
\beq{Rval2}
R_x^{\ell_0+1}(y) = O\bfrac{\deg^-(u)}{m (np)^2}.
\eeq

{\bf (ii) $e$ is a back edge inducing a directed cycle.}

In this case, there is an edge $e=(u,v)$ where $u \in Y_i, v \in Y_j$ and $j >i$.
Let $vPu$ denote the path from $v$ to $u$ in $\Tu_y$, and $C$ the cycle $vPu, (u,v)$.
%Repeating the discussion used for Type 1 walks in $Q_x^{\ell_0+1}(y)$,
There is some $k\geq 1$ such that the walk is
$P_0=xPu,(uCu)^k,uPy$. Let $\s$ be the length of $C$,
let $\t$
be the distance from $u$ to $y$ in $\Tu_y$, and let $s=\t+k\s$.
Let $\ell=\ell_0-s$. Then $\ell+1$ is the length of the walk $xPu$ from $x$ to $u$ prior to
the final $s$ steps.

Either $\ell < (1+4\eta/5)\La$ and the
 argument in Case 1 ($i\geq 4\eta/5)\La$) above can be applied, giving us the bound
\beq{gggg}
R_x^{\ell_0+1}(y)= O\bfrac{(100/\eta)^{\ell_3}{\D_0}^{9\eta\La/10}}{(c_0np)^{\ell_0+1}}=
O(n^{-1-\eta/20}).
\eeq
Or $\ell \ge (1+4\eta/5)\La$ and
we adapt Case 2. Let $w$ be the predecessor of $u$ on $P_0$.
We can use Remark \ref{rem2}
as above to obtain $P_x^{\ell}(w) \le\ooi \deg^-(w)/m$.
As $k\s \ge 2, \t \ge 0$, (the worst case is $u=y, w \in N^-(y)$), we obtain
\beq{Rval3}
R_x^{\ell_0+1}(y)= O \bfrac{\deg^-(w)}{m(np)^2}.
\eeq
Thus, using \eqref{Rval1}, \eqref{Rval2}, \eqref{gggg}, \eqref{Rval3} we have
\beq{Rval}
R_x^{\ell_0+1}(y)= O\brac{ \frac{1}{n}\cdot \frac{1}{(np)^2}}.
\eeq

We have therefore shown that
$S_x^{\ell_0+1}(y)+Q_x^{\ell_0+1}(y)+R_x^{\ell_0+1}(y)=o(1/n)$ completing the proof that
\beq{anot}
P_x^{\ell_0+1}(y)\leq(1+o(1))\frac{\deg^-(y)}{m}
\eeq

\blem{kl1}
For all $y\in V$,
$$\p_y=(1+o(1))\frac{\deg^-(y)}{m}.$$
\elem
\proofstart
It follows from \eqref{lowpi} that for any $y\in V$,
\begin{equation}\label{finally2}
\p_y=\sum_{x\in V}\p_xP^{(\ell_0+1)}_x(y)\leq
(1+o(1))\frac{\deg^-(y)}{m}\sum_{x\in
  V}\p_x=(1+o(1))\frac{\deg^-(y)}{m}.
\end{equation}
The lemma now follows from Lemma \ref{kl}.
\proofend

\section{Stationary distribution: Removing Assumption 1}\label{remove}
\subsection{Large average degree case}\label{bigdegree}
\subsubsection{$np\geq n^\d$.}

We can deal with this case by using a concentration
inequality \eqref{KimVu} from
Kim and Vu \cite{KV}: Let $\Upsilon=(W,E)$
be a hypergraph where $e\in E$ implies that $|e|\leq s$. Let
$$Z=\sum_{e\in E}w_e\prod_{i\in e}z_i$$
where the $w_e,e\in E$ are positive reals and the $z_i,i\in W$ are
independent random variables taking values in $[0,1]$.
For $A\subseteq W,|A|\leq s$ let
$$Z_A=\sum_{\substack{e\in E\\e\supseteq A}}w_e\prod_{i\in e\setminus A}z_i.$$
Let $M_A=\E(Z_A)$ and $M_j(Z)=\max_{A, |A|\geq j}M_A$ for $j\ge0$.
There exist positive constants $a$ and $b$
such that for any $\l>0$,
\begin{equation}\label{KimVu}
\Pr(|Z-\E(Z)|\geq a\l^s\sqrt{M_0M_1})\leq b|W|^{s-1} e^{-\l}.
\end{equation}
For us, $W$ will be the set of edges of $\vec{K}_n$
the complete digraph on $n$ vertices. {$z_i$ will be
the indicator variable for the presence of the $i$th edge
of $\vec{K}_n$.} $E$ will be the set of
sets of edges in walks of length $s=\rdup{2/\d}$ between two fixed vertices $x$ and
$y$ in $\vec{K}_n$, and $w_e=1$. $Z$ will be the number of
walks of length $s$ that are in $D_{n,p}$. In which case we have
\begin{align*}
&\E(Z)=(1+o(1))n^{s-1}p^s\\
&M_j\leq
(1+o(1))n^{s-j-1}p^{s-j}\le(1+o(1))\E(Z)/np\qquad for\ j\geq 1.
\end{align*}
So $M_0=\E(Z)$ and
applying \eqref{KimVu} with $\l=(\log n)^2$ we see that for any
$x,y$ we have
$$\Pr(|Z-\E(Z)|=O(\E(Z)n^{-\d/2}\log^{O(1)} n))= 1-O(n^{-3}).$$
Thus \whp\
$$P_x^s(y)=(1+o(1))\frac{n^{s-1}p^s}{((1-\e_1)np)^s}\sim
\frac{1}{n},\qquad\forall x,y\in V.$$
We now finish with the arguments of Lemmas \ref{kl} and \ref{kl1}.
\proofend

\subsection{Small average degree case}\label{Smavdeg}
\subsubsection{Lower bound on stationary distribution}

A vertex is
{\em small} if it has in-degree or out-degree at most $np/20$ and
{\em large} otherwise.
In the proofs of Section \ref{lbd} we assumed $x,y$ were large.
We proceed as in Section \ref{lowerbound} but initially restrict our
analysis to large $x,y$. Also, with the exception of $Y_1$ we do not
include small vertices when creating
the $X_i,Y_i$.
Avoiding the $\leq n^{1/5}$ small
vertices (see Lemma \ref{newlem0}(a)) is easily incorporated because in the proof we have allowed
for the avoidance of $n^{2/3+o(1)}$ vertices from $\bigcup_iX_i$ etc.
Provided there are no small vertices in $N^-(y)$, our previous lower bound analysis holds.
In this way, we show  for all large $x,y$ that,
\beq{cucu}
P_x^{(2\ell+1)}(y)\geq (1-o(1))\frac{\deg^-(y)}{m}.
\eeq
If $x$ is small, then it will only have large out-neighbours (see Lemma \ref{newlem0}(c)) and
so if $y$ is large then
\beq{moo}
P_x^{(2\ell+2)}(y)=\frac{1}{\deg^+(x)}\sum_{z\in
  N^+(x)}P_z^{(2\ell+1)}(y)\geq (1-o(1))\frac{\deg^-(y)}{m}.
\eeq
A similar argument deals with small $y$ and $x$ arbitrary i.e.
\beq{mooo}
P_x^{(2\ell+2)}(y)=\sum_{z\in
  N^-(y)}\frac{P_x^{(2\ell+1)}(z)}{\deg^+(z)}\geq (1-o(1))\sum_{z\in
  N^-(y)}\frac{\deg^-(z)}{m}\frac{1}{\deg^+(z)}\geq (1-o(1))\frac{\deg^-(y)}{m}.
\eeq
We have used Lemma \ref{newlem0}(e) to justify the last inequality.

\ignore{
We thus obtain a proof that for all  $x,y$,
$$P_x^{(2\ell+2)}(y)\geq (1-o(1))\frac{\deg^-(y)}{m}.$$
(We apply the argument of \eqref{mooo} to \eqref{cucu} to replace $\ell+1$ by $\ell+2$).
}

In the case that some $u\in N^-(y)$ has small out-degree,
then by Lemma \ref{newlem0}(c) there is at most one
such $u$ \whp. For $z \in N^-(y)$, we repeat the argument above for each factor $P_x^{2\ell+1}(z)$.
The extra term $\varsigma^*(y)$ now arises from $\deg^-(u)/\deg^+(u)$ and
$$P_x^{(2\ell+2)}(y)=\sum_{z\in N^-(y)}\frac{P_x^{(2\ell+1)}(z)}{\deg^+(z)}\geq (1-o(1))\frac{1}{m}
\sum_{z\in N^-(y)}\frac{\deg^-(z)}{\deg^+(z)}\geq (1-o(1))\frac{\deg^-(y)+\varsigma^*(y)}{m}.$$
We can  now proceed as in \eqref{finally}.

\subsubsection{Upper bound on stationary distribution}
We first explain how  the upper bound proof  in Section \ref{upperbound} alters if Assumption 1
is removed.
The assumption that the minimum degree was at least $c_0np$ was used in the following places:
\begin{enumerate}
\item We assumed in Section \ref{lbd} that $\deg^+(x), \deg^-(y) \ge c_0np$.
These assumptions can be circumvented by using
Lemma \ref{newlem0}(c) with the methods used in the lower bound case.

\ignore{
\item In \eqref{c_0npval}, \eqref{smallcase1}.  If we assume that $y$ is large and
has no small in-neighbours, then the proof of Section \ref{ubd} is unaltered.
}

\item In \eqref{Sval}, \eqref{Pit}, \eqref{Rval1}. In these cases
 we used $(c_0np)^{\ell_0}$ as a lower bound on the product of out-degrees on a path of length
$\l$ for some $\l\geq \ell_1$. Using Lemmas \ref{newlem0} and \ref{newlem2},  we see that
 small vertices are at weak distance at least $\ell_{10}$ and therefore there can be at most 11
such vertices on any walk length $\ell_0+1$.
Thus, after dropping Assumption 1, we replace this lower bound by $(c_0np)^{\l-11}$, and
 the proof  continues essentially unchanged.

\item In the proof of Lemma \ref{xlem2} we made a re-scaling $B=1000/(c_0np)^{2\ell+1}$.
The exponent $2 \ell +1$ was replaced by $\ell_0+1$ in the
proof of \eqref{Zconcx} in Lemma \ref{xlem5}.
We now replace $\ell_0+1$ by $\ell_0-10$.

\item In the proof of Lemma \ref{yLem} we made a re-scaling $U(i)=W(y,i)\cdot(c_0np)^{i}$
at each level $3 \le i \le \ell$.
Assume that $2\ell_2<\ell_{10}$ i.e. $\eta\leq 1/250$
so that there is at most one small vertex  $u$ in $Y$.
If we replace
$(c_0np)^i$ by
$(c_0np)^{i-1}$ does not affect our concentration results, provided $i \ge 3$.
The bounds on $U_v$ are now $(c_0/np)(c_0/C_0)^i \le U_v \le 1$, and $\e_i= 1/\sqrt{(A \log n)^{i-2}}$.
If the small vertex $u \in N^-(y)$ then the direct calculations
used in the lower bound hold.  If the small vertex $u$ is in levels $i=2,3$ this
adds an extra term of $O(\deg^-(u)/(m (np)^{i-1}))$ to our estimate of $Z_x^{\ell_0+1}(y)$
in  Section \ref{upperbound}.

\item In \eqref{PType1}, \eqref{Rval2}, \eqref{Rval3}.
It follows from Lemma \ref{newlem2}, that if e.g.
$\Tu_y$ contains a non-tree edge, then no vertex of $\Tu_y$
is small, and the calculations in the proof are
unaltered.
\end{enumerate}

Thus
the proof as is works perfectly well if we assume that $y$ is large and if it
has no small in-neighbours and there is no small vertex
in $Y$.
We call such  a vertex $y$ {\em ordinary}.

If $y$ is small then from Lemmas \ref{newlem0} and \ref{newlem2} we can
assume that all of its in-neighbours are ordinary.
This is under the assumption that $2\ell_2<\ell_{10}$ e.g. if $\eta\leq 1/250$.
So in this case we can use
Lemma \ref{newlem0}(e) and obtain
$$P_x^{(\ell_0+2)}(y)=\sum_{\xi\in N^-(y)}\frac{P_x^{(\ell_0+1)}(\xi)}{\deg^+(\xi)}
\leq\frac{1+o(1)}{m} \sum_{\xi\in
N^-(y)}\frac{\deg^-(\xi)}{\deg^+(\xi)}=(1+o(1))\frac{\deg^-(y)}{m}.$$

Suppose now that $y$ is large and that there is a
small vertex $u\in Y$. We can assume from Lemma \ref{newlem2} that $Y$
does not contain any edge not in $\Tu_y$. Either
$u \in N^-(y)$ or, if not, from point 4. of the discussion above,
an extra $O(\deg^-(u)/(m (np)))$ is added to $Z_x^{\ell_0+1}(y)$
for the probability of the $(x,y)$-walk going via $u$.

In the case where $u \in N^-(y)$ then as in the lower bound
\begin{align*}
P_x^{(\ell_0+1)}(y)&\le \frac{1+o(1)}{m}\brac{\frac{\deg^-(u)}{\deg^+(u)}+ \sum_{u\in N^-(y)\setminus
w}\frac{\deg^-(u)}{\deg^+(u)}}\\
&= \frac{(1+o(1))}{m}\brac{\deg^-(y)+\varsigma^*(y)}.
\end{align*}
We have now
completed the proof of the asymptotic steady state without Assumption 1.

\section{Mixing time and the conditions of Lemma \ref{MainLemma}}
\label{smallsec}
\subsection{Upper Bound on Mixing time}\label{secmix}
Let $T$ be a mixing time as defined in \eqref{4} and
let $\ell=O(\log_{np} n)$ be given by \eqref{ell}.
We prove that (\whp) $T$
 satisfies
\begin{equation}\label{T}
T=o(\ell\log n)=o((\log n)^2).
\end{equation}

Define
\begin{equation}\label{hdbar}
\hd(t)=\max_{x,x'\in V}|P_{x}^{(t)}-P_{x'}^{(t)}|
\end{equation}
to be the maximum over $x,x'$ of the variation distance between $P_{x}^{(t)}$ and $P_{x'}^{(t)}$.
It is proved in Lemma 20 of Chapter 2 of Aldous and Fill \cite{AlFi}
that
\beq{afil}
\hd(s+t)\leq \hd(s)\hd(t)\mbox{ and }\max_x|P_{x}^{(t)}-\p_x|\leq \hd(t).
\eeq

Equation \eqref{qs3} implies that {whp}
\begin{equation}\label{Qxz}
\hd(2\ell+1)=o(1),
\end{equation}
and so \eqref{T} follows immediately from  \eqref{afil} and \eqref{Qxz}.

\subsection{Conditions of Lemma \ref{MainLemma}}
We see immediately from \eqref{T} that
Condition (b) of Lemma \ref{MainLemma} is satisfied.

We  show below that \whp\ for all $v\in V$
\begin{equation}\label{RT}
R_T(1)=1+o(1).
\end{equation}

Using \eqref{RT}, the proof that  Condition (a) of  Lemma  \ref{MainLemma} is  satisfied,
is as follows.
Let $\l=1/KT$ as in \eqref{lamby}.
%The value of $T=o(\log^2 n)$ is given by \eqref{T}.
For $|z|\leq 1+\l$,  we have
$$R_T(z)\geq 1-\sum_{t=1}^Tr_t|z|^t\geq 1-(1+\l)^T\sum_{t=1}^Tr_t=1-o(1).$$
Thus  for $v\in V$, the value of $p_v$ in \eqref{pv} is given by
\beq{pvval}
p_v=(1+o(1)) \frac{\deg^-(v)}{m}.
\eeq

{\bf Proof of \eqref{RT}:}
If $d\geq (\log n)^2$, then the minimum out-degree of $D_{n,p}$ is
$\Omega(d\log n)$. In which case we have for any $x,y$
\begin{equation}\label{uprob}
\Pr(\cW_v(t)=y\mid \cW_v(t-1)=x)=O\brac{\frac{1}{d\log n}}.
\end{equation}
The expected number of returns to $v\in V$ by $\cW_v$ during $T$ steps,
is therefore $O(T/d\log n){=o(1)}$.

Now assume that $d\leq (\log n)^2$.
\begin{enumerate}[(i)]
\item Lemma \ref{smallsets} implies that if $H$ is the subgraph of $D_{n,p}$
induced by vertices at weak distance at most $\Lambda/20$ from $v$ then $H$ contains
at most $|V(H)|$ edges.
\item Lemma \ref{newlem0} implies that there is at most one small vertex in $H$.
\item Lemma \ref{newlem2} implies that there is no
small vertex within weak distance 10 of a
weak cycle of length $\leq 10$.
\end{enumerate}
Assume that conditions (i), (ii), (iii) hold.
Let $A_4$ denote the set of vertices $u\neq v$
such that $D_{n,p}$ has a path of length at most 4 from $u$ to $v$.
We show next that:
\beq{yo3}
\text{With probability }1-O(1/(np)^2),\ \cW_v(i)\notin A_4,\,1\leq i\leq 4.
\eeq
For this to happen, there has to be a cycle $C$ of
length at most 8 containing $v$. If such a cycle exists then
all vertices within weak distance 10 of $v$ have degree
at least $np/20$. Furthermore, the only way that the walk
can reach $A_4$ in 4 or less steps is via this cycle.
This verifies \eqref{yo3}. Assume then that $\cW_v(i)\notin A_4,\,1\leq i\leq 4$.

Suppose next that there is a time $T_1\leq T$ such that
$\cW_v(T_1)=v$. Let\\ $T_2=\min\set{\t\leq T_1:\cW_v(t)\in A_4,\t\leq t\leq T_1}$.
It must be the case that $d(T_2)=4$ where $d(t)$ is the distance from $\cW_v(t)$ to $v$.

If $A_4$ does not contain a small weak cycle then
the walk must proceed directly to $v$ in 4 steps. The
probability of this is $O(1/(np)^3)$, since at most one vertex on the path of length
4 from $x=\cW(T_2)$ to $v$
will be of degree at most $np/20$.

If there is a small weak cycle $C$ then there is an edge $e$ of $C$ whose
removal leaves an in-branching
of depth 4 into $v$. There are now 2 paths that $\cW$ can follow
from $x$ to $v$. One uses $e$ and one does not. Each path has a probability
of $O(1/(np)^4)$ of being followed. Putting this altogether we see
that the expected number of returns to $v$ is $O(1/(np)^2+T/(np)^3)=o(1)$.
This completes the proof of \eqref{RT}.

\section{The Cover Time of $D_{n,p}$}\label{cta}
\subsection{Upper Bound on the Cover Time}
For $np=d \log n$, $d $ constant, let
$t_0=(1+\e)\brac{d\log\bfrac{d}{d-1}}n\log n$.
For $np=d\log n$ $d=d(n) \rai$ let
$t_0=(1+\e)n\log n$. In both cases we assume $\e \to 0$
sufficiently slowly to ensure that all inequalities  below are valid.

Let $T_D(u)$ be the time taken by the random walk $\cW_u$
to visit every vertex of $D$. Let $U_t$ be
the number of vertices of $D$ which have not been visited by
$\cW_u$ at step $t$.
We note the following:
\begin{eqnarray}
\gap{.4}C_u=\E(T_D(u))&=& \sum_{t > 0} \Pr(T_D(u) \ge t), \label{ETG} \\
\label{TGa}
\Pr(T_D(u)\geq t)=\Pr(T_D(u) > t-1)&=&\Pr(U_{t-1}>0)\le \min\{1,\E(U_{t-1})\}.
\end{eqnarray}
Recall that
 $\ul A_v(t)$ denotes the event that $\cW_u(t)$ did
not visit vertex $v$ in the interval $[T,t]$.
It follows from (\ref{ETG}), (\ref{TGa}) that for any $t\geq T$,
\begin{equation}\label{shed1}
C_u \le t+1+ \sum_{s \ge t} \E(U_{s})\le t+1+\sum_v\sum_{s \ge t}\Pr(\ul A_v(s)).
\end{equation}

Assume first that $d(n)\to\infty$.
If  ${s/T\to\infty}$ then \eqref{frat} of Lemma \ref{MainLemma}
{together with the value of $p_v$
given by \eqref{pvval}, and concentration of in-degrees }implies that
\begin{equation}\label{Prv}
\Pr(\ul A_v(s))\leq
(1+o(1))\exp\set{-\frac{(1-o(1))s}{n}}+O(e^{-\OM(s/T)}).
\end{equation}
Plugging \eqref{Prv} into \eqref{shed1} we get
\begin{eqnarray}
C_u&\leq&t_0+1+2n\sum_{s\geq
  t_0}\brac{\exp\set{-\frac{(1-o(1))s}{n}}+O(e^{-\OM(s/T)})}\label{placeof}\\
&\leq&t_0+1+3n^2\exp\set{-\frac{(1-o(1))t_0}{n}}+O(nTe^{-\OM(t_0/T)})\nn\\
&=&(1+o(1))t_0.\nn
\end{eqnarray}

We now assume that $d$ is {bounded as $n\to\infty$}, and the conditions of Lemma
 \ref{vertexdegrees} hold.
For $v\in V$ we
have
$$\Pr(\ul A_v(s))=(1+o(1))\exp\set{-(1+o(1/\log n))\p_vs}+O(e^{-\OM(s/T)})$$
where, {by Lemma \ref{kl}},
$$\p_v\geq (1-o(1))\frac{\deg^-(v)}{m}.$$
In place of \eqref{placeof} we  use the bounds on the number of vertices
of degree $k$ given in Lemma
\ref{vertexdegrees}, {in terms of the sets $K_i, \; i=0,1,2,3$}. Thus
\beq{Cuq}
C_u\leq t_0+1+o(1)+\sum_{i=0}^3S_i
\eeq
where
\begin{eqnarray*}
S_i&=&\sum_{k\in K_i}D(k)\sum_{s\geq t_0}\exp\set{-\frac{(1-o(1))ks}{m}}\\
&\leq&2m\sum_{k\in K_i}\frac{D(k)}{k}e^{-(1-o(1))kt_0/m}\\
&\leq&2m\sum_{k\in K_i}\frac{D(k)}{k}\bfrac{d-1}{d}^{(1+\e/2)k}.
\end{eqnarray*}
The main term {occurs at $i=3$.  Using \eqref{aproxval}, \eqref{K3val}},
{the fact that $(nep(d-1))/(kd))^k$ is maximized
at $k=np(d-1)/d$, and $m=dn \log n \ooi$ \whp, we see that}
\begin{eqnarray}
S_3&\leq&\frac{8m}{n^{d-1}}\sum_{k=c_0np}^{\D_0}
\bfrac{nep}{k}^k\bfrac{d-1}{d}^{(1+\e/2)k}\nonumber\\
&\leq&8m \; \D_0 \;
e^{-\e c_0np/2d}\nonumber\\
&=&o(t_0).
\end{eqnarray}

{Note that $K_0=0$.
We next consider the cases $i=1,2$.
For $i=1$, we refer first to  Lemma \ref{vertexdegrees}(i-a).
If  $d-1\geq (\log n)^{-1/3}$ then $K_1=\emptyset$.
If $d-1< (\log n)^{-1/3}$, then
 $D(k)\leq (\log\log n)^2$, from \eqref{K1val}.
In this case $t_0=O((1/(d-1))) dn \log n$.
 Thus
\begin{eqnarray}
S_1&\leq&m\sum_{k\in K_1}\frac{D(k)}{k}\bfrac{d-1}{d}^{(1+\e/2)k}\nonumber\\
&\leq&m\sum_{k=1}^{15}\frac{(\log\log n)^2}{k}\bfrac{d-1}{d}^{(1+\e/2)k}\nonumber\\
&=& O(t_0)(\log\log n)^2 (d-1)^{-\e/2} \nonumber \\
&=&o(t_0)\label{S1}
\end{eqnarray}

For $i=2$, by Lemma \ref{vertexdegrees} if $d-1< (\log n)^{-1/3}$ and $k \ge 16$,
and using \eqref{K2val} we have $D(k)\leq (\log n)^4$.
 Thus
\begin{eqnarray}
S_2&\leq&m\sum_{k\in K_2}\frac{D(k)}{k}\bfrac{d-1}{d}^{(1+\e/2)k}\nonumber\\
%&\leq&m\sum_{k\in K_2}\frac{(\log n)^4}{k}\bfrac{d-1}{d}^{(1+\e/2)k}\nonumber\\
&\le& O(t_0)\sum_{k\in K_2}\frac{\log^4 n}{k}(d-1)\bfrac{d-1}{d}^{(1+\e/2)k}\nonumber\\
&=& O(t_0) \log^4 n (\log n)^{-(19/3+\e/8)}\nonumber\\
&=&o(t_0). \label{S2.1}
\end{eqnarray}
If
 $d-1\geq (\log n)^{-1/3}$ then  by Lemma \ref{vertexdegrees}(i-a) $\min \{k \in K_2\} \ge
(\log n)^{1/2}$, and $|K_2|=O(\log\log n)$. Thus, as $d$ is bounded
\begin{eqnarray}
S_2
&=& O(t_0)\sum_{k\ge (\log n)^{1/2} }\frac{\log \log n}{k}(d-1)\bfrac{d-1}{d}^{(1+\e/2)k}\nonumber\\
&=&o(t_0)\label{S2}
\end{eqnarray}

The upper bound on cover time of  $C_u\leq t_0+o(t_0)$ now follows from \eqref{Cuq}--\eqref{S2}.
}
\subsection{Lower Bound on the Cover Time}

For $np=d \log n$, let
$t_1=(1-\e)d\log\bfrac{d}{d-1}n\log n$.
Here $\e\to 0$
sufficiently slowly so that all inequalities claimed  below are valid.

{\bf Case 1: $np\leq n^\d$ where $0<\d\ll\eta$ is a positive constant.}\\
Let $k^*=(d-1)\log n$, and let $V^{*}=\set{v:\;\deg^-(v)=k^*\ and\ \deg^+(v)=d\log n}$.
\Whp the size $|V^{*}|\ge n^*=\frac{n^{\g_d}}{4\p\log n(d(d-1))^{1/2}}$ (see Lemma
\ref{vertexdegrees}(ii)).
Let us first work assuming $d\geq 1.05$.
In this case $\g_d=(d-1)\ln(d/(d-1))\geq .15$ and we write $n^*=n^{\g_d-o(1)}$.
The maximum degree in $D$ is at most $\D_0=O(np)$ and so
$V^{*}$ contains a
sub-set $V_1^{*}$ of size $n^{\g_d/2}$ such that $v,w\in V_1^{*}$ and $x\in V$
implies
\begin{eqnarray}
dist(x,v)+dist(x,w)&>& {\Lambda/100}.\label{PP1}\\
dist(y,x)+dist(x,y)&>& {\Lambda/50},\ for\ y=v,w.\label{PP2}
\end{eqnarray}
Here "dist" refers to directed distance in $D_{n,p}$ and recall that $\Lambda=\log_{np}n$.

 Each $v\in V_1^*$ has $\p_v\sim \frac{d-1}{dn}$
and so we can choose a subset $V^{**}$ of size $\geq n^{\g_d/3}$ such that if
$v_1,v_2\in V^{**}$ then
\begin{equation}\label{close}
|\p_{v_1}-\p_{v_2}|\leq \frac{1}{n\log^{10} n}.
\end{equation}
Indeed, suppose that $\p_v\in \left[\frac{d-1}{2dn},\frac{2(d-1)}{dn}\right]$ for $v\in V_1^*$.
Divide this interval into $\log^{10} n$ equal sized sub-intervals and
then use the pigeon-hole principle.

\ignore{
In addition, if $d\leq \log^3n$ then we can assume that if $v\in V^{**}$ then there
is no edge joining two vertices
within out-distance $\leq (\log\log n)^2$ of $v$, see Lemma \ref{lemcov}.
}
Now choose $u\notin
V^{**}$ and let $V^{\dagger}$ denote
the set of vertices in $V^{**}$ that have not been visited by
$\cW_u$ by time
$t_1$.
Then $\E(|V^{\dagger}|) \rai$, as the following calculation shows;
\[
\E(|V^{\dagger}|)\geq n^{\g_d/3}\brac{\exp\set{-\frac{(1+o(1))k^*t_1}{m}}-o(e^{-\OM(t_1/T)})}-T,
\]
where the last term accounts for possible visits before time $T$.

Now assume that $1+o(1)\leq d\leq 1.05$. In these circumstances we have $n^*=\log^\om n$ where
$\om\to\infty$, see \eqref{EVL}.
Equations \eqref{PP1}, \eqref{PP2} now hold for all $v,w\in V^*$. This follows
from Lemma \ref{newlem0} because the vertices of $V^*$ are small.
The size of $V^{**}$ is at least $n^*/(\log n)^{10}$
and we can again write
\begin{eqnarray*}
\E(|V^{\dagger}|)&\geq&
\frac{n^*}{(\log n)^{10}}\brac{\exp\set{-\frac{(1+o(1))k^*t_1}{m}}-o(e^{-\OM(t_1/T)})}-T\\
&\to&\infty.
\end{eqnarray*}
As in previous papers, {see for example \cite{CFreg}, we will
finish our proof by using}, the Chebyshev inequality to
show that $V^{\dagger}\neq \emptyset$ \whp, thus completing the proof of Theorem
\ref{thmain}.
This will follow if we can prove that
$$\var(|V^\dagger|)=o(\E(|V^\dagger|^2)+O(|V^{**}|^2n^{-2})=o(\E(|V^\dagger|)^2).$$
To establish this inequality, we will show that if $v,w\in V^{**}$ then
\begin{equation}\label{need}
\Pr(\ul A_v(t_1)\cap\ul A_w(t_1))\leq (1+o(1))\Pr(\ul A_v(t_1))\Pr(\ul A_w(t_1)).
\end{equation}
To prove this, we identify vertices $v,w$ into a
``supernode'' $\s$ to obtain a digraph $D_\s$ with $n-1$ vertices. In
this digraph $\s$ has in-degree $deg^-(v)+deg^-(w)=2k^*$ {and out-degree $2d \log n$}.

{\bf The stationary distribution of $D_\s$.}\\
Let $\p^*$ denote the vector of steady states in
$D_\s$.
The arguments we used in Sections \ref{structure} and \ref{SteadyState}
remain valid in $D_{\s}$, and thus
$$\p^*_\s\sim (1-o(1))\frac{2k^*}{m}.$$
However, we need to be  more precise.
For a vertex $x$ of $D_\s$ let
$$\hp_x=\begin{cases}\p_x&x\neq\s\\\p_v+\p_w&x=\s\end{cases}.$$
We will prove  for all $x\in V(D_\s)$, that
\begin{equation}\label{pdiff}
|\p^*_x-\hp_x|=O\bfrac{1}{n(\log n)^{8}}.
\end{equation}
{\bf Proof of \eqref{pdiff}.}\\
Let $\xi=\hp-\p^*$ be the difference between  $\hp$ and $\p^*$.
Let $P^*$ be the transition matrix of the walk on $D_\s$, then
$$P^*(x,y)=\begin{cases}P(x,y)
&x,y\neq\s\\(P(v,y)+P(w,y))/2&x=\s\\P(x,v)+P(x,w)&y=\s\end{cases}.$$
Let $\xi'$ be the transpose of $\xi$.
It follows from the steady state equations that
$$(\xi'P^*)_x=\begin{cases}\hp_x-\p^*_x&x\notin N^+(\{v,w\})\\
\hp_x-\p^*_x+\frac{\p_w-\p_v}{2}P(v,x)&x\in N^+(v)\\
\hp_x-\p^*_x+\frac{\p_v-\p_w}{2}P(w,x)&x\in N^+(w)\end{cases}.$$
We rewrite this as
\beq{home1}
\xi'(I-P^*)=\eta'
\eeq
where $\eta_x=0$ for $x\notin N^+(\{v,w\})$ and $|\eta_x|\leq |\p_v-\p_w|/2$ otherwise.

Multiplying \eqref{home1} on the right by $M=\sum_{t=0}^{T-1}(P^*)^t$ we have
\begin{equation}\label{closer}
\xi'(I-P^*)M=\xi'(I-(P^*)^T)=\eta'M.
\end{equation}
{Let}
\beq{closerx}
(P^*)^T=\Pi+E
\eeq
where $\Pi$ is the $(n-1)\times (n-1)$ matrix with each row equal to $(\p^*)'$.
The definition of $T$ implies that each entry of $E$ has
absolute value bounded by $n^{-3}$.

Now write $\xi=\a\p^*+\z$ where $\z\perp\p^*$. It follows from
$(\p^*)'P^*=(\p^*)'$ and \eqref{closer}
that
$$(\a \p^*+\z)'(I-(P^*)^T)=\z'(I-(P^*)^T)=\z'(I-\Pi-E)=\eta'M.$$

Now
$$\z'(I-E)=\z'(I-(P^*)^T+\Pi)=\eta'M+\z'\Pi.$$
As $\z\perp\p^*$ this implies that
\beq{fff}
\z'(I-E)\z=\eta'M\z.
\eeq
Note that
\beq{ffff}
|\eta'M\z|\leq \sum_{t=0}^{T-1}|\eta'(P^*)^t\z|\leq T|\eta||\z|,
\eeq
{where $|z|$ denotes the $\ell_2$ norm of $z$.}

Now
\beq{fffff}
|\z'(I-E)\z|\geq |\z|^2-|\z'E\z|\geq |\z|^2-n^{-3}
\brac{\sum_{i=1}^{n-1}|\z_i|}^2\geq |\z|^2(1-n^{-2}).
\eeq
It follows from \eqref{fff}, \eqref{ffff} and \eqref{fffff} that
$$
|\z|^2(1-n^{-2})\leq T|\eta||\z|
$$
and so using \eqref{close} we find that
\beq{ggg}
|\z|=O\bfrac{1}{n(\log n)^{8}}.
\eeq
Now let {\bf 1} denote the $(n-1)$-vector of 1's. Then
$$0=1-1=(\hp-\p^*)'{\bf 1}=\xi'{\bf 1}=\a+\z'{\bf 1}.$$
Using \eqref{ggg} this gives
$$
|\a|\leq  |{\bf 1}|\,|\z|=O\bfrac{1}{n^{1/2}(\log n)^{8}}.
$$
Now
$\xi_x=\a\p^*_x+\z_x$ for all $x$ and so
$$\xi_x^2\leq 2\a^2(\p^*_x)^2+2\z_x^2=
O\brac{\frac{1}{n(\log n)^{16}}\cdot \frac{1}{n^2}
+\frac{1}{n^2(\log n)^{16}}}=O\bfrac{1}{n^2(\log n)^{16}}.$$
This completes the proof of \eqref{pdiff}. \proofend

{\bf Proof of \eqref{need}.}\\
For  $v\in V^{**}$, we first tighten \eqref{RT}  to
\begin{equation}\label{RTT}
R_v=1+o(1/(\log n)^2).
\end{equation}
Assume first that $np\leq \log^{10}n$. Then
\eqref{PP1} and \eqref{PP2} imply that for $1\leq t\leq (\log n)^{2/3}$, vertex $v$
will be at distance $\geq 2\log^{2/3} n-t$ from $\cW_v(t)$. Then once the walk is
at a vertex $w$ within distance $\log^{2/3} n$ of $v$
its chance of getting closer is only $O(1/\log n)$.
This being true with at most one exception for a vertex of low out-degree.
The probability that there is a time $t$ such that $\cW_v$ is within $\log^{2/3} n$ of $v$
and it makes 10 steps closer to $v$ in the next 100 steps
is $O(T/\log^9 n)=O(1/\log^7 n)$. This implies
\eqref{RTT}. If $np\geq \log^{10}n$ then we use $R_v\leq 1+(1+o(1))T/np$.

Similarly,
\begin{equation}\label{RTTs}
R_\s=1+o(1/(\log n)^2).
\end{equation}
The mixing time $T$ in what follows is the maximum of the mixing times for $D$ and the
maximum over $v,w$ for $D_\s$.
Using the suffix $\Pr_\s$ to denote probabilities related to random walks in $D_\s$
and using \eqref{pdiff}, it follows that
\begin{eqnarray}
\Pr_\s(\ul A_\s(t_1))&\leq&\exp\set{-\frac{(1+O(T\p^*_\s))\p_{\s}^*t_1}{m}}-o(e^{-\OM(t_1/T)})
\nonumber\\
&\leq&\exp\set{-\frac{(1+o(1/\log n))(\p_{v}+\p_w)t_1}{m}}-o(e^{-\OM(t_1/T)})\nonumber\\
&=&(1+o(1))\Pr(\ul A_v(t_1))\Pr(\ul A_w(t_1)).\label{omega1}
\end{eqnarray}
But, using rapid mixing in $D_\s$,
\begin{eqnarray*}
\Pr_\s(\ul A_\s(t_1))&=&\sum_{x\neq \s}P^T_{\s,u}(x)\Pr_\s(\cW_x(t)\neq \s,1\leq t\leq t_1-T)\\
&=&\sum_{x\neq \s}(\p^*_x+O(n^{-3}))\Pr_\s(\cW_x(t)\neq \s,1\leq t\leq t_1-T)
\end{eqnarray*}
On the other hand,
\begin{eqnarray*}
\Pr(\ul A_v(t_1)\cap \ul A_w(t_1))
&=&\sum_{x\neq v,w}P^T_u(x)\Pr(\cW_x(t)\neq v,w,T\leq t\leq t_1)\\
&=&\sum_{x\neq v,w}(\p_x+O(n^{-3}))\Pr(\cW_x(t)\neq v,w,1\leq t\leq t_1-T)
\end{eqnarray*}
But,
$$\Pr_\s(\cW_x(t)\neq \s,T\leq 1\leq t_1-T)=\Pr(\cW_x(t)\neq v,w,1\leq t\leq t_1-T)$$
because random walks from $x$ that do not meet $v,w$ or $\s$ have the same measure in both digraphs.

It follows that
\begin{align}
&\Pr(\ul A_v(t_1)\cap \ul A_w(t_1))-\Pr_\s(\ul A_\s(t_1))\nonumber\\
&=\sum_{x\neq v,w}(\p_x-\p^*_x+O(n^{-3}))\Pr(\cW_x(t)\neq v,w,1\leq t\leq t_1-T)\nonumber\\
&\leq O\bfrac{1}{n\log^8n}\sum_{x\neq v,w}\Pr(\cW_x(t)\neq v,w,1\leq t\leq t_1-T)\nonumber\\
&\leq O\bfrac{1}{n\log^8n}\sum_{x\neq v,w}\frac{P^T_u(x)}{P^T_u(x)}
\Pr(\cW_x(t)\neq v,w,1\leq t\leq t_1-T)\nonumber\\
&\leq O\bfrac{1}{n\log^8n}O(n\log n)\sum_{x\neq v,w}P^T_u(x)
\Pr(\cW_x(t)\neq v,w,1\leq t\leq t_1-T)\nonumber\\
\noalign{since $P^T_u(x)=\Omega(1/n\log n)$}\nonumber\\
&\leq O\bfrac{1}{\log^7n}\Pr(\ul A_v(t_1)\cap \ul A_w(t_1)).\label{omega2}
\end{align}

Equations \eqref{omega1} and \eqref{omega2} together imply \eqref{need}.
\proofend

{\bf Case 2: $np\geq n^\d$.}\\
In this range we take $t_1=(1-\e)n\log n$ and let $V^{*}$ be the set of vertices of degree
$\rdown{np}$. A simple second moment calculation shows that \whp\ we
have $|V^{*}|=\Omega((np)^{1/2-o(1)})$. We then choose $\e$ so that
$\E(|V^\dagger|)\geq (np)^{1/4}$. It is then only a matter of verifying
\eqref{need}. The details are as in the previous case.

This completes the proof of Theorem \ref{thmain}.
\proofend

{\bf Acknowledgement:} We thank several referees whose
insight and hard work has helped to make this paper correct and (hopefully)
more readable,

\end{document}